\documentclass[11pt]{article}
\usepackage{amsmath,amsfonts,amsthm}
\usepackage{geometry}
\usepackage{amssymb}
\usepackage{graphics}
\usepackage{graphicx}

\newtheorem{prethm}{Spectral Convergence:  Theorem}
\newtheorem{prethm2}{Theorem}
\newtheorem{bigthm}[prethm]{Heat Kernel Convergence:  Theorem}
\newtheorem{bigthm2}[prethm2]{Theorem}
\newtheorem{defglap}{Definition}
\newtheorem{defac}[defglap]{Definition}
\newtheorem{defcon}[defglap]{Definition}
\newtheorem{resblow}[defglap]{Definition}
\newtheorem{accsingle}[defglap]{Definition}
\newtheorem{accon}[defglap]{Definition}

\newtheorem{aclemma1}{Lemma}

\newtheorem{bmap}[defglap]{Definition}
\newtheorem{fibration}[defglap]{Definition}
\newtheorem{defbmani}[defglap]{Definition}
\newtheorem{defbheat}[defglap]{Definition}

\newtheorem{bheatcomp}{Technical Theorem}
\newtheorem{heatcalc}[defglap]{Definition}
\newtheorem{composition}[bheatcomp]{Technical Theorem}

\newtheorem{acheatcalc}[defglap]{Definition}
\newtheorem{accomposition}[prethm]{Technical Theorem}
\newtheorem{acparametrix}[prethm]{Technical Theorem}

\newtheorem{accheatcalc}[defglap]{Definition}

\newcommand{\cS}{\mathcal{S}}
\newcommand{\cG}{\mathcal{G}}
\newcommand{\cC}{\mathcal{C}}
\newcommand{\cL}{\mathcal{L}}

\newcommand{\cH}{\mathcal{H}}
\newcommand{\cT}{\mathcal{T}}

\newcommand{\N}{\mathbb{N}}

\newcommand{\R}{\mathbb{R}}

\newcommand{\heart}{\begin{flushright} $\heartsuit$ \end{flushright}}

\begin{document}

\title{Spectral Geometry and Asymptotically Conic Convergence}
\author{Julie Marie Rowlett \thanks{Email:  rowlett@math.ucsb.edu.  Department of Mathematics, South Hall Room 6607, University of California, Santa Barbara, California 93106.  AMS Classification primary:  58J50, secondary:  58J35.}}
\maketitle

\begin{abstract}

In this paper we define \em asymptotically conic convergence \em in which a family of smooth Riemannian metrics degenerates to have an isolated conic singularity.  This convergence is related to, and generalizes, the analytic surgery metric degeneration of \cite{hmm}.  For a conic metric $(M_0, g_0)$ and a scattering metric $(Z, g_z)$ we define a nonstandard blowup, the \em resolution blowup, \em in which the conic singularity in $M_0$ is resolved by $Z.$  Equivalently, the resolution blowup resolves the boundary of the scattering metric using the conic metric; the resolution space is a smooth compact manifold.  This blowup induces a smooth family of metrics $\{ g_{\epsilon} \}$ on the compact resolution space $M,$ and we say $(M, g_{\epsilon})$ converges \em asymptotically conically \em to $(M_0 , g_0)$ as $\epsilon \to 0.$   

Let $\Delta_{\epsilon}$ and $\Delta_0$ be geometric Laplacians on $(M, g_{\epsilon})$ and $(M_0, g_0),$ respectively. Our first result is convergence of the spectrum of $\Delta_{\epsilon}$ to the spectrum of $\Delta_0$ as $\epsilon \to 0.$  Note that this result implies spectral convergence for the $k$-form Laplacian under certain geometric hypotheses.  This theorem is proven using rescaling arguments, standard elliptic techniques, and the $b$-calculus of \cite{tapsit}.  Our second result is technical:  we construct a parameter ($\epsilon$) dependent heat operator calculus which contains, and hence describes precisely, the heat kernel for $\Delta_{\epsilon}$ as $\epsilon \to 0.$  The consequences of this result include:  the existence of a polyhomogeneous asymptotic expansion for $H_{\epsilon}$ as $\epsilon \to 0,$ with uniform convergence for all time, and a precise description of the behaviour of $H_{\epsilon}$ uniformly as $\epsilon \to 0.$  To prove this result we construct several manifolds with corners (heat spaces) using both standard and non-standard blowups, on which we construct corresponding heat operator calculi.  A parametrix construction modeled after the heat kernel construction of \cite{tapsit} and a maximum principle type argument complete this proof.  

\end{abstract}

\section{Introduction}
Cheeger and Colding wrote a series of three papers between 1997 and 2000 on the Gromov-Hausdorff limits of families of smooth, connected Riemannian manifolds with lower Ricci curvature bounds \cite{cc1}, \cite{cc2}, \cite{cc3}.  They proved Fukaya's conjecture of 1987 \cite{fukaya}:  on any pointed Gromov-Hausdorff limit space  of a  family $\{M^n _i \}$ of connected Riemannian manifolds with Ricci curvature bounded below, a self adjoint extension of the scalar Laplacian can be defined, with discrete spectrum in the compact case, so that the eigenvalues and eigenfunctions of the scalar Laplacians on $M^n _i$ behave continuously as $i \to \infty.$  An example, provided by Perelman, showed that the results of Cheeger and Colding could not be extended to the $k$-form Laplacian or to more general geometric Laplacians.  In 2002, Ding proved convergence of the heat kernels and Green's functions in the same setting \cite{ding}.  The estimates are uniform for time bounded strictly away from zero. These results are impressive;  the only hypothesis is a lower Ricci curvature bound.  It would be useful to prove a more general spectral convergence result for geometric Laplacians and to obtain uniform estimates on the heat kernels for all time.  In order to obtain such results, it becomes necessary to impose more structure on the manifolds and the way in which they converge to a singular limit space.  

Let $M$ be a fixed compact manifold with Riemannian metric and let $H$ be an embedded orientable hypersurface with defining function $x$ and smooth metric $g_H,$ and let
$$g_{\epsilon} := dx^2 +  (\epsilon^2 + x^2 ) g_H \quad \epsilon \in [0, 1).$$
As $\epsilon \to 0$, $g_{\epsilon} \to x^2 g_H + dx^2,$ which has an isolated conic singularity at $x=0$.  Geometrically, $M$ is pinched along the hypersurface $H$ as $\epsilon \to 0,$ and the resulting metric has a conic singularity as $x \to 0$.  The study of this metric collapse is the content of the $1990$ thesis of McDonald \cite{pat}.

In 1995, Mazzeo and Melrose developed pseudodifferential techniques to describe the behavior of the spectral geometry under another specific type of metric collapse known as \em analytic surgery \em similar to the conic collapse considered by McDonald.  As $\epsilon \to 0$ the metrics
$$g_{\epsilon} = \frac{|dx|^2}{x^2 + \epsilon^2} + h \, \to \, \frac{|dx|^2}{x^2} + h = g_0;$$
$g_0$ is an exact $b$-metric on the compact manifold with boundary $\bar{M}$ obtained by cutting $M$ along $H$ and compactifying as a manifold with boundary, hence the name, \em analytic surgery\em.  Under certain assumptions on the associated Dirac operators,\footnote{In a later collaboration with Hassel \cite{hmm} these hypotheses were removed.} Mazzeo and Melrose proved
$$ \lim_{\epsilon \to 0} \eta (\partial_{\epsilon} ) = \eta_b (\partial_{\bar{M}} ),$$
where $\eta_b (\partial_{\bar{M}} )$ is the $b$-version of the eta invariant introduced by Melrose.  These results were proven by analyzing the resolvent family of the Dirac operators $\partial_{\epsilon}$ uniformly near zero.  This led to a precise description of the behaviour of the small eigenvalues.  Our goal is to obtain precise analytic results like those of \cite{hmm}.

The convergence considered here, \em asymptotically conic (ac) convergence, \em is more restrictive than that of Cheeger, Colding and Ding, but more general than the conic degeneration of McDonald.   We note, however, that ac convergence does not require Ricci curvature bounds. The conic collapse of McDonald, the analogous smooth collapse of a higher codimension submanifold and the collapse of an open neighborhood of the manifold with some restrictions on the local geometry all fit this new definition.  Before stating our spectral convergence results for geometric Laplacians, we recall their definition.
\begin{defglap}
Let $(E, \nabla)$ be a Hermitian vector bundle over a Riemannian manifold $(M, g)$ with metric-compatible connection $\nabla.$  A geometric Laplacian is an operator $\Delta$ acting on sections of $E$ which has the form
$$\Delta = \nabla^* \nabla + R,$$
where $R$ is a non-negative self-adjoint endomorphism of $E.$  By the Weitzenb\"ock Theorem, the Laplacian on $k$-forms is a geometric Laplacian, as is the Hodge Laplacian and the conformal Laplacian; any geometric Laplace-type operator is a geometric Laplacian.
\end{defglap}

Our results are the following.

\begin{prethm}
\label{prethm}
Let $(M_0, g_0 )$ be a compact Riemannian $n$-manifold with isolated conic singularity,  and let $(Z, g_z)$ be an asymptotically conic space, with $n \geq 3.$  Assume $(M, g_{\epsilon} )$ converges asymptotically conically to $(M_0 , g_0).$  Let $(E_0, \nabla_0)$ and $(E_z, \nabla_z)$ be Hermitian vector bundles over $(M_0, g_0)$ and $(Z, g_z),$ respectively, so that each of these bundles in a neighborhood of the boundary is the pullback from a bundle over the cross section $(Y, h).$  Let $\Delta_0,$ $\Delta_z$ be the corresponding Friedrich's extensions of geometric Laplacians, and let $\Delta_{\epsilon}$ be the induced geometric Laplacian on $(M, g_{\epsilon}).$  Assume $\Delta_z$ has no $\mathcal{L}^{2}$ nullspace.  Then the accumulation points of the spectrum of $\Delta_{\epsilon}$ as $\epsilon \to 0$ are precisely the points of the spectrum of $\Delta_0,$ counting multiplicity.  
\end{prethm}

The setting for our next result is the acc heat space, a manifold with corners constructed in section 7.  

\begin{bigthm}
Let $(M_0, g_0 )$ be a compact Riemannian $n$-manifold with isolated conic singularity,  and let $(Z, g_z)$ be an asymptotically conic space, with $n \geq 3.$  Assume $(M, g_{\epsilon} )$ converges asymptotically conically to $(M_0 , g_0).$  Let $(E_0, \nabla_0)$ and $(E_z, \nabla_z)$ be Hermitian vector bundles over $(M_0, g_0)$ and $(Z, g_z),$ respectively, so that each of these bundles in a neighborhood of the boundary is the pullback from a bundle over the cross section $(Y, h).$  Let $\Delta_0,$  $\Delta_z$ be the corresponding Friedrich's extensions of geometric Laplacians, and let $\Delta_{\epsilon}$ be the induced geometric Laplacian on $(M, g_{\epsilon}).$  Then the associated heat kernels $H_{\epsilon}$ have a full polyhomogeneous expansion as $\epsilon \to 0$ on the asymptotically conic convergence (acc) heat space with the following leading terms at the $\epsilon =0$ faces:
$$H_{\epsilon} (z, z', t) \to H_0 (z, z', t')   \textrm{ at } F_{0101} $$
$$H_{\epsilon} (z, z', t) \to (\rho_{1010})^2 H_b (z, z', \tau)  \textrm{ at } F_{1010},$$
$$H_{\epsilon} (z, z', t) = O\left( (\rho_{1001})^2 (\rho_{0110})^2 \right) \textrm{ at } F_{1001}, \, F_{0110}.$$

Above, $\rho_{wxyz}$ is the defining function for boundary face $F_{wxyz},$ $H_0$ is the heat kernel for $\Delta_0,$ $H_b$ is the $b$-heat kernel for $\Delta_z$ and $t'$ and $\tau$ are rescaled time variables with 
$$t' = \frac{t}{(\rho_{1010})^2}, \, \tau = \frac{t}{(\rho_{1001} \rho_{0110})^2}.$$   This convergence is uniform in $\epsilon$ for all time and moreover, the error term is bounded by $C_N \epsilon t^N,$ for any $N \in \N_0.$    

\end{bigthm}

\textbf{Remark:}  The convergence of the heat kernels $H_{\epsilon} (z, z', t)$ on $M \times M \times \R^+$ to the heat kernel $H_0 (z, z', t')$ on $M_0$ follows immediately from Theorem 2.  The error term is $\epsilon^2 (H_b (z, z', \tau) + O(1) ) + O(\epsilon t^{\infty})$ as $\epsilon, t \to 0$ and $\epsilon^2 (H_b (z, z', \tau) + O(1)) + O(\epsilon)$ as $\epsilon \to 0$ for all $t.$  However, Theorem 2 is stronger than this statement:   the theorem specifies the leading order behavior of $H_{\epsilon}$ precisely at the $\epsilon =0$ faces of the acc heat space and gives the existence of a full polyhomogeneous expansion as $\epsilon \to 0.$  

These theorems are proven in sections 5 and 7, respectively.  In section 2 we define the resolution blowup and ac convergence.  Sections 3 and 4 contain a brief review of geometric and analytic results and terminology on manifolds with corners.  In section 6 we construct the heat spaces and heat operator calculi that will be used to prove the main theorem in section 7.    

This work is based on the author's doctoral dissertation completed at Stanford University in June, 2006 under the supervision of Rafe Mazzeo.  The author wishes to thank Rafe Mazzeo for excellent advising, Andras Vasy for many helpful conversations and suggestions, and Richard Melrose for insightful comments.  

\section{Asymptotically Conic Convergence}
The definition of \em asymptotically conic convergence \em involves three geometries:  a family of smooth metrics on a compact manifold, a conic metric on a compact (incomplete) manifold, and an \em asymptotically conic  \em or ac scattering metric on a compact (complete) manifold with boundary.  

First, we define $ac$ scattering metric.
\begin{defac}
An $ac$ scattering metric $(\bar{Z}, g_z)$ is a smooth metric on a compact $n$-manifold with boundary, $\partial \bar{Z} = (Y, h)$ a smooth compact $(n-1)$ manifold.   $\bar{Z}$ has a product decomposition $(0, r_1)_r \times Y$ in a neighborhood of the boundary defined by $r=0,$ so that on this neighborhood, 
$$g_z = \frac{dr^2}{r^4} + \frac{h(r)}{r^2}, \quad h(r) \to h \textrm{ as } r \to 0,$$
in other words, $h$ extends to a $\mathcal{C}^{\infty}$ tensor on $[0, r_1).$
Uniquely associated to the ac scattering metric $(\bar{Z}, g_z)$ is the complete non-compact manifold $Z$ with asymptotically conic end; $Z$ is known as an \em ac space.\em\footnote{Note that asymptotically conic spaces are sometimes called ``asymptotically locally Euclidean," or $ALE$.  However, that term is often used for the more restrictive class of spaces that are asymptotic at infinity to a cone over a quotient of the sphere by a finite group, so to avoid confusion, we use the term \em asymptotically conic. \em}  Letting $\rho = \frac{1}{r},$ there is a compact subset $K_z \subset Z$ so that 
$$Z - K_z \cong (1/r_1, \infty)_{\rho} \times Y, \quad  g_z |_{Z - K_z} = d\rho^2 + \rho^2 h(1/\rho).$$
\end{defac}
A familiar example of an ac scattering metric is the standard metric on the radial compactification of $\R^n$ with  boundary $\mathbb{S}^{n-1}$ at infinity.  

Next, we define the conic metric on a compact manifold. 
\begin{defcon}
Let $M_c$ be a compact metric space with Riemannian metric $g$.  Then, $(M_c, g)$ has an isolated conic singularity at the point $p$ and $g$ is called a \em conic metric \em if the following hold.
\begin{enumerate}
\item $(M_c - \{p\} , g )$ is a smooth, open manifold. 
\item There is a neighborhood $N$ of $p$ and a function $x : N - \{p\} \to (0,x_1]$ for some $x_1 > 0,$ such that $N - \{p\}$ is diffeomorphic to $(0, x_1]_x \times Y$ with metric $g = dx^2 + x^2 h(x)$ where $(Y, h)$ is a compact, smooth $n-1$ manifold and $\{ h(x) \}$ is a smooth family of metrics on $Y$ converging to $h$ as $x \to 0,$ in other words, $h$ extends to a $\mathcal{C}^{\infty}$ tensor on $[0,x_1)_x \times Y.$   
\end{enumerate}
\end{defcon}
Associated to a manifold with isolated conic singularity is the manifold with boundary  obtained by blowing up the cone point, adding a copy of $(Y, h)$ at this point.  Then, $x$ is a boundary defining function for the boundary, $(Y, h).$  We use $M_0 ^0 = M_c - \{p\}$ to denote the smooth incomplete conic manifold, $M_c$ to denote the metric space closure of $M_0 ^0,$ and $M_0$ to denote the associated manifold with boundary.  We work almost exclusively with $M_0.$  

The \em resolution blowup \em is a nonstandard blowup that resolves an isolated conic singularity in a compact manifold using an asymptotically conic space.  In the definition of resolution blowup we use the notation $M \cup_{\phi} N$ for a smooth manifold constructed from the smooth manifolds $M$ and $N$ with a diffeomorphism $\phi$ from $V \subset N$ to $U \subset M$ that gives the equivalence relation, $V \ni p \sim \phi(p) \in U.$  $M \cup_{\phi} N$ is the disjoint union of $M$ and $N$ modulo the equivalence relation of  $\phi.$  The smooth structure on $M \cup_{\phi} N$ and the topology is induced by that of $M$ and $N.$  

\begin{resblow} \label{resblow}
Let $(M_0, g_0)$ be a compact $n$ manifold with isolated conic singularity and let $(Z, g_z)$ be an asymptotically conic space of dimension $n,$ so that $\partial M_0 = \partial \bar{Z} = (Y,h).$  Then, 
$$M_0 ^0 = K_0 \cup V_0,$$
where $V_0 \cong (0,x_1)_x \times Y,$ and $K_0$ is compact. With this diffeomorphism 
$$g_0 = dx^2 + x^2 \tilde{h}(x) \quad \textrm{ on $(0, x_1) \times Y.$}$$
We may assume the boundary of $K_0$ in $M_0 ^0$ is of the form $\partial K_0 = \{ x = x_1 \} \cong Y,$ and we may extend $x$ smoothly to $K_0$ so that $2 x_1 > x \geq x_1$ on $K_0.$   Similarly, 
$$Z = K_z \cup V_z,$$
where $V_z \cong (\rho_1, \infty)_{\rho} \times Y,$ and $K_z$ is compact.  With this diffeomorphism
$$g_z = d\rho^2 + \rho^2 h(\rho , y) \quad \textrm{ on $(\rho_1, \infty) \times Y.$}$$
We may assume that $\partial K_z$ is of the form $\partial K_z = \{ \rho = \rho _1 \} \cong Y,$ and we similarly extend $\rho$ smoothly to $K_z$ so that $\rho_1 \geq \rho > \rho_1 /2$ on $K_z.$   

Let $\delta = \textrm{ min } \{x_1, 1/\rho_1 \}.$  Then for $0< \epsilon < \delta,$ and $R > \rho_1,$ let 
$$M_{0, \epsilon} =   \{(x, y) \in M_0 :  x > \epsilon \} , \quad Z_R =  \{ (\rho, y) \in Z : \rho < R \} .$$  
The  \em \textbf{resolution blowup} \em of $(M_0 , g_0)$ by $(Z, g_z)$ is, 
$$M_{\epsilon} := M_{0, \epsilon} \cup_{\phi_{\epsilon}} Z_{1/\epsilon},$$
where the joining map $\phi$ is defined for each $\epsilon$ by
$$\phi_{\epsilon} : M_{0,\epsilon} - M_{0, \delta} \to Z_{1/\epsilon} - Z_{1/\delta} , \quad \phi_{\epsilon} (x, y) = \left( \frac{ x}{\delta \epsilon} , y \right).$$
For $\delta > \epsilon > \epsilon' >0,$ the manifolds $M_{\epsilon}$ and $M_{\epsilon'}$ are diffeomorphic, and so the resolution blowup of $M_0$ by $Z$ which we call $M$ is unique up to diffeomorphism.  
\end{resblow}

\textbf{Remark:}  The resolution blowup resolves the singularity in a conic manifold using an ac scattering space.  Instead of resolving the singularity in $M_0$ using $Z,$ we may equivalently define the resolution blowup to resolve the boundary of $\bar{Z}$ using $M_0.$\footnote{Let $r = 1/ \rho$ be the defining function for $\partial \bar{Z}.$  The resolution blowup of $\bar{Z}$ by $M_0$ is
$$M_0 \cup_{\psi_{\epsilon}} \bar{Z}, \quad \psi_{\epsilon} (x, y) = \left( \frac{\epsilon \delta}{x} , y \right),$$
where the patching map $\psi_{\epsilon}$ is defined on $M_0 - K_0$ with image $(0, \delta)_r \times Y$ a neighborhood of $\partial \bar{Z}.$  The resulting smooth compact resolution space is diffeomorphic to $M.$}

The ac single space, analogous to the analytic surgery single space in \cite{eta}, is the setting for the definition of ac convergence.  
 
\begin{accsingle} \label{single}
Let $(M_0, g_0)$ be a conic metric and let $(\bar{Z}, g_z)$ be a scattering metric; assume both are dimension $n$ with the same cross section $(Y, h)$ at the boundary, and assume $\delta=1$ (definition \ref{resblow}).  Then, $M_0 \cong  ( (0,1)_x \times Y ) \cup K_0$ and $\bar{Z} \cong ( (0,1)_r \times Y) \cup K_z$ with $\partial K_0 \cong Y \cong \partial K_z.$  The \em asymptotically conic convergence (acc) single space $\cS$ \em  is defined as follows,
$$\mathcal{S} := [0,1)_x \times [0,1)_r \times Y \cup (K_0 \times \{ x =1, r \neq 1 \} ) \cup (K_z \times \{r=1, x \neq 1 \}).$$
\end{accsingle}

The smooth structure of $\cS$ is induced by that of $M_0$ and $\bar{Z}.$  Namely, smooth functions on $\cS$ are functions which are smooth jointly in $x$ and $r$ on $(0,1)_x \times (0,1)_r \times Y,$ smoothly extend to a smooth function on $K_0$ at $x=1,$ on $K_z$ at $r =1,$ on $M_0$ at $r=0,$ and on $\bar{Z}$ at $x=0.$   To give a precise description of the metrics considered here, we define below the acc tensor, a smooth, polyhomogeneous, symmetric 2-cotensor on the acc single space.  

\begin{accon} \label{accon}
Let $\mathcal{S}$ be the acc single space associated to $M_0$ and $\bar{Z}$ as in definition \ref{single}.  Let $\epsilon(p) = x(p) r(p): \cS \to [0,1),$ where $x, r$ are extended to $K_0$ and $K_z$ respectively to be identically $1.$  We define the \em acc tensor \em $\mathcal{G}$ as follows:
\begin{displaymath}
\mathcal{G} = \left \{ \begin{array}{ll}  dx^2 + x^2 h(x) + (xr)^2 \left( \frac{dr^2}{r^4} + \frac{h(r)}{r^2} \right)    & x, r \in [0,1) \\
dx^2 + x^2 h(x) + (xr)^2(  g_z |_{K_z} )   & r=1, \, x \in [0,1) \\
 g_0 |_{K_0} + (xr)^2 \left( \frac{dr^2}{r^4} + \frac{h(r)}{r^2} \right)  & x=1, \,  r \in [0,1) 
\\
\end{array} \right .
\end{displaymath}

For $0< \epsilon < 1,$ let $M_{\epsilon} = \{ xr = \epsilon \} \subset \mathcal{S};$ note that this $M_{\epsilon}$ is diffeomorphic to the resolution blowup $M$ of $M_0$ by $Z.$\footnote{In the definition of the acc single space we have assumed $\delta$ from the definition of resolution blowup is $1.$  This is to simplify our calculations.  On $\{ xr = \epsilon \}$ we then have $r = \frac{\epsilon}{x},$ equivalently $x = \frac{\epsilon}{r},$ so letting $\delta =1$ in the definition of resolution blowup the identification of $\{ xr = \epsilon\} \subset \cS$ and the resolution blowup $M_{\epsilon}$ follows immediately.}  The family of metrics $\{ g_{\epsilon} = \cG |_{M_{\epsilon}} \}$ on $M$ is said to converge \em asymptotically conically \em to $(M_0, g_0).$  
\end{accon}

\textbf{Remarks}
\begin{enumerate}
\item The acc single space has two boundary hypersurfaces at $\epsilon =0;$ these are diffeomorphic to $M_0$ at $r=0$ and $\bar{Z}$ at $x=0,$ and they meet in a codimension $2$ corner diffeomorphic to $Y.$  There are also boundary hypersurfaces at $\epsilon =1$ which we ignore since we are interested in $\epsilon \to 0.$ 

\item Since $g_0$ is a smooth metric on $M_0 \cong \left( (0,1)_x \times Y \right) \cup K_0 $ and $g_z$ is a smooth metric on $\bar{Z} \cong  \left( (0,1)_r \times Y  \right) \cup K_z,$ $\cG$ is a smooth symmetric $2$-cotensor on $\cS$ which extends smoothly across neighborhoods where it is piecewise defined.  

\item At $r=0,$ $\cG$ restricts to $\cG |_{ \{ r=0 \} } = g_0;$ as $x \to 0,$ $\cG$ vanishes to order $2.$

\item On $M_{\epsilon} \subset \cS$ where $0 < r(p), x(p) < 1,$ 
$$r = \frac{\epsilon}{x} \, \quad \implies dr^2 = \frac{\epsilon^2}{x^4} dx^2,$$
so
$$g_{\epsilon} = dx^2 + x^2 h(x) = \epsilon^2 \left(\frac{dr^2}{r^4} + \frac{h(r)}{r^2} \right).$$

\item On $M_{\epsilon}$ for $x \equiv \epsilon,$ $r \equiv 1$ and this subset is diffeomorphic to $K_z$ with metric $g_{\epsilon} = \epsilon^2 g_z |_{K_z}.$  On $M_{\epsilon}$ for $x \equiv 1,$ $r \equiv \epsilon$ and this subset is diffeomorphic to $K_0$ with metric $g_{\epsilon} = g_0 + O(\epsilon)^2.$  

\end{enumerate}

The following lemma is useful both for visualizing ac convergence and for proving spectral convergence.

\begin{aclemma1} \label{aclemma1}
Let $(M_0, g_0)$ and $(Z, g_z)$ be as in definitions \ref{resblow}, \ref{single}, \ref{accon} and let $(M, g_{\epsilon})$ converge asymptotically conically to $(M_0, g_0).$  Then, there exists a family of diffeomorphisms $\{ \phi_{\epsilon} \}$ from a fixed open proper subset $U \subset M$ to increasing neighborhoods $Z_{1/\epsilon} \subset Z$ such that $g_{\epsilon} |_U \cong \left( \epsilon^2 (\phi_{\epsilon})^* g|_{Z_{1/\epsilon}} \right)|_U. $   Moreover, on $M-U,$ $g_{\epsilon} \to g_0$ smoothly as $\epsilon \to 0$ and any $K \subset \subset M_0 ^0$ is diffeomorphic to some fixed $K' \subset M$ so that $g_{\epsilon} \to g_0$ smoothly and uniformly on $K'.$  
\end{aclemma1}

\textbf{Proof}  

The existence of $\phi_{\epsilon}$ and $U \subset M$ follows immediately from the definition of resolution blowup and the diffeomorphism between the resolution blowup $M$ and $\{ xr = \epsilon \} \subset \cS.$  By the above remarks, on the neighborhood $U \subset M$ where this diffeomorphism is defined, $g_{\epsilon} |_U = (\epsilon^2) (\phi_{\epsilon})^* g_z |_{Z_{1/\epsilon}}.$

Since $g_{\epsilon} = g_0 + \mathcal{O}(\epsilon^2),$ the smooth convergence of $g_{\epsilon}$ to $g_0$ on $M-U$ follows immediately.  Any compact subset $K \subset \subset M_0$ is contained in $M_{0, \epsilon}$ for some $\epsilon >0$ and so is diffeomorphic to $K_{\epsilon} \subset \subset M_{\epsilon}$ and to $K' \subset \subset (M - U).$  Conversely, any $K \subset \subset (M-U)$ is diffeomorphic to $K_{\epsilon} \subset M_{\epsilon}$ and to $K' \subset \subset M_{0, \epsilon} \subset M_0.$    
 \heart

\section{Geometric Preliminaries}

This section is a brief review of the theory and terminology of manifolds with corners, $b$ maps, and blowups.  A complete reference is \cite{tapsit}, see also \cite{edge}.  

\subsection{Manifolds with Corners}
Let $X$ be a manifold with corners.  This means that near any of its points, $X$ is modeled on a product $[0, \infty)^k \times \R^{n-k}$ where $k$ depends on the point and is the maximal codimension of the boundary face containing that point.  We also assume that all boundary faces of $X$ are embedded so they too are manifolds with corners.  The space $\mathcal{V} (X)$ of all smooth vector fields on $X$ is a Lie algebra under the standard bracket operation.  It contains the Lie subalgebra
\begin{equation} \label{btangent} \mathcal{V}_b (X) := \{ V \in \mathcal{V}(X) \, ; \, V \textrm{ is tangent to each boundary face of $X.$} \} \end{equation}
Then $\mathcal{V}_b (X)$ is itself the space of all smooth sections of a vector bundle, 
$$\mathcal{V}_b (X) = \mathcal{C}^{\infty} (X; \, ^bTX),$$
where $^b TX$ is the bundle defined so that the above holds and is called the $b$-tangent bundle.  

\subsubsection{Blowing Up}
An embedded codimension $k$ submanifold $Y$ of a manifold with corners $X$ is called a $p$-submanifold ($p$ for product) if near each point of $Y$ there are local product coordinates so that $Y$ is defined by the vanishing of some subset of them.  In other words, $X$ and $Y$ must have consistent local product decompositions.  Then one can define a new manifold with corners $[X; Y]$ to be the normal blowup of $X$ around $Y.$  This is obtained by replacing $Y$ by its inward-pointing spherical normal bundle.  The union of this normal bundle and $X - Y$ has a unique minimal differential structure as a manifold with corners so that the lifts of smooth functions on $X$ and polar coordinates around $Y$ are smooth.  One can also consider iterated blowups, written $\left[[X; Y_1] ; Y_2]\right],$ where $Y_1$ is a $p$-submanifold of $X$ and $Y_2$ is a $p$-submanifold of $Y_1.$  We may consider any finite sequence of such blowups.  If we have such a sequence of embedded $p$-submanifolds,
$$X \supset Y_1 \supset Y_2 \supset Y_3 \ldots \supset Y_n$$
then the iterated blowup 
$$\left[ [X; Y_1]; Y_2];  \ldots ; Y_n \right]$$
can be performed in any order with the same result \cite{eta}.  Blowups may also be defined using equivalence classes of curves \cite{tapsit}.  Let $r$ be a defining function for the $p$ submanifold and consider the family of curves $\gamma(t) = (r(t), y(t))$ such that
$$\gamma(t) \in Y \iff t=0,$$
$$r(t) = O(t).$$
Let $E$ be the set of equivalence classes of all such curves with 
$$\gamma \sim \gamma' \iff (y-y')(t) = O(t), \, (r-r')(t) = O(t^2).$$
There is a natural $\R^+$ action on $E$ given by
$$\R^+ \ni a: \gamma(t) \to \gamma(at).$$
Then $E$ modulo this equivalence relation is naturally diffeomorphic to $N^+(Y),$ the inward pointing spherical normal bundle of $Y,$ so we can define $[X; Y]$ by
$$[X; Y] = (X-Y) \cup E/(\R^+ - \{0\}).$$
We can also define \em parabolic \em blowups in certain contexts \cite{epmm}.  Let $Y$ be a $p$-submanifold of codimension $k$ so that there exist local coordinates $(r, y) = (r_1, \ldots, r_k, y_1, \ldots, y_{n-k})$ in a neighborhood of $Y$ with $r_i$ vanishing precisely at $Y,$ and so that $dr_1$ induces a sub-budle of the tangent bundle $TX.$  Instead of the above equivalence classes of curves we consider $\gamma(t)$ such that
$$\gamma(t) = (r_1(t), \ldots, r_k(t), y_1(t), \ldots y_{n-k} (t)) \in Y \iff t=0,$$
$$r_i(t) = O(t), \, i \neq 1, \, r_1(t) = O(t^2).$$
Two such curves are equivalent if
$$\gamma \sim \gamma' \iff (y_j-y_j')(t) = O(t), (r_i-r_i')(t) = O(t^2) \, i \neq 1, \, (r_1 (t) - r_1 ' (t) ) = O(t^3).$$
Since $dr_1$ is a sub-bundle of $TX$ there is a natural $\R^+$ action on the set of equivalence classes $E_2$ of all such curves,
$$\R^+ \ni a : \gamma(t) \to (r_1( a^2 t), r_i( at), \ldots, y_j (at)).$$
The set of equivalence classes of all such curves modulo this $\R^+$ action is naturally diffeomorphic to the inward pointing $r_1$-parabolic normal bundle of $Y,$
$$E_2 / (\R^+ - \{0 \}) \cong PN^+ _{r_1} (Y).$$
We define the $r_1$-parabolic blowup of $X$ around $Y$ as the union of $X-Y$ and this inward pointing $r_1$-parabolic bundle,  
$$[X; Y, dr_1] := (X-Y) \cup PN^+ _{r_1} (Y).$$  
The union of this $r_1$-parabolic bundle and $X-Y$ again has a unique minimal differential structure as a manifold with corners so that the lifts of smooth functions on $X$ and $r_1$-parabolic coordinates around $Y$ are smooth.  By $r_1$-parabolic coordinates around $Y,$ we mean the coordinates,
$$\rho = (r_1 ^2 + r_2 ^4 + \ldots + r_k ^4)^{\frac{1}{4}}, \,  \theta = (\theta_1, \ldots \theta_n) \in \mathbb{S}^{n-1},$$
with local coordinates $(r_1, \ldots r_k, y_1, \ldots, y_{n-k})$ in a neighborhood of $Y$ satisfying
$$r_i = \rho \theta_i, \, i \neq 1, \, r_1 = \rho^2 \theta_1,  \quad y_j = \rho \theta_j.$$

For any parabolic or spherical blowup there is a natural blow-down map $\beta_*: [X; Y] \to X$ and corresponding blow-up map $\beta^* : X \to [X; Y],$ so that the image of $Y$ under $\beta^*$ is a boundary hypersurface of $[X; Y]$ diffeomorphic to the inward pointing spherical (or parabolic) normal bundle of $Y.$  As such, 
$$[X; Y] = (X-Y) \cup \beta^* (Y).$$

\subsubsection{$b$-Maps and $b$-Fibrations}
\begin{bmap} 
Let $M_1$ be a manifold with boundary hypersurfaces, $\{ N_j\}_{j=1} ^k ,$ and defining functions $r_j .$  Let $M_2 $ be a manifold with boundary hypersurfaces, $\{ L_i \}_{i=1} ^l ,$ and defining functions $\rho_i .$  Then $f: M_1 \to M_2$ is called a \em $b$-map \em if for every $i$ there exist nonnegative integers $e(i, j)$ and a smooth nonvanishing function $h$ such that $f^* (\rho_i) = h \prod_{j=1} ^k r_j ^{e(i, j) } .$  
\end{bmap}

The image under a $b$-map of the interior of each boundary hypersurface of $M_1$ is either contained 
in or disjoint from each boundary hypersurface of $M_2 $ and the order of vanishing of the differential of $f$ is constant along each boundary hypersurface of $M_1 .$  The matrix $(e(i, j))$ is called the lifting matrix for $f.$  

In order for the map, $f,$ to preserve polyhomogeneity, stronger conditions are required.  Associated to a manifold with corners are the $b$-tangent and cotangent bundles, $^b TM$ (\ref{btangent}) and $^b T^* M.$\footnote{These are also called the totally characteristic tangent and cotangent bundles.}  The map $f$ may be extended to induce the map $^b f_* : ^b TM_1 \to ^b TM_2 .$  

\begin{fibration}
The $b$-map, $f : M_1 \to M_2 ,$ is called a \em $b$-fibration \em if the associated maps $^b f_*$ at each $p \in \partial M_1$ are surjective at each $p \in \partial M_1$ and the lifting matrix $(e(i, j))$ has the property that for each $j$ there is at most one $i$ such that $(e(i, j)) \neq 0.$  In other words, $f$ does not map any boundary hypersurface of $M_1$ to a corner of $M_2 .$  
\end{fibration}

\subsubsection{$b$-manifolds and the $b$-blowup}

A $b$-manifold is a manifold with corners that is closely related to conic manifolds and ac scattering metrics. 
\begin{defbmani} \label{bmani}
Let $(X , g)$ be a smooth Riemannian manifold with boundary $(Y, h)$ and boundary defining function $x$ such that in a collared neighborhood $N$ of the boundary $X$ has a product decomposition, $N \cong [0,x_1)_x \times Y$ and in this neighborhood 
$$g = \frac{dx^2}{x^2} + h(x),$$
where $h(x)$ is a smoothly varying family of metrics on $Y$ that converges smoothly to $h$ as $x \to 0.$  Then $(X, g)$ is said to be a \em $b$-manifold.  \em  
\end{defbmani}

Equivalently, a $b$-manifold is a complete manifold with asymptotically cylindrical ends.  The Schwartz kernels of operators on a $b$-manifold with reasonable regularity lift to a blown up manifold called the $b$-double space.  This space is obtained from $X^2$ by performing a radial blowup called the $b$-blowup along the codimension $2$ corner at the boundary in each copy of $X$ and it is written $X^2 _b,$
\begin{equation} \label{bdouble} X^2 _b = [X \times X; \partial X \times \partial X] = [X \times X ; Y \times Y].\end{equation}

For any manifold $M$ with boundary having a product structure in a neighborhood of the boundary, we may define the $b$-blowup in the analogous way, $M^2 _b := [M \times M ; \partial M \times \partial M].$

\subsection{Asymptotically conic convergence double space}

The acc double space is an instructive model for the more sophisticated acc heat space construction in section 7.  Let 
$$\mathcal{S}^2 _b := [\mathcal{S}^2; C_1 \times C_1],$$
where we recall $C_1 \cong Y$ is the codimension two corner in each copy of $\cS.$  The acc double space $\mathcal{D}$ is the submanifold of $\mathcal{S}^2 _b$ defined by the vanishing set of $f(p) = x(p)r(p) - x'(p) r'(p) = \epsilon(p) - \epsilon ' (p):$
$$\mathcal{D} = \{ p \in \mathcal{S}^2 _b: f(p) = x(p) r(p) - x'(p) r'(p) = 0\} = \{ p \in \cS^2 _b : \epsilon(p) = \epsilon'(p) \}.$$
The acc double space has various boundary faces but we are only interested in those at $\epsilon =0.$  There are four boundary faces at $\epsilon =0,$ described in the following table.  Here and throughout, we label each face $F_{wxyz}$ where the subscript indicates the order to which each of the scalar variables $x, r, x', r'$ vanishes at that face.   

\begin{tabular}{|l|l|l|}
\hline
Arising from & face & geometry \\
\hline
$x=0, x'=0,$ & $F_{1010} $ & $[ \bar{Z} \times \bar{Z} ; Y \times Y]$ \\
$x=0, r'=0$ & $F_{1001} $ & $[ \bar{Z} \times M_0; Y \times Y]$ \\
$r=0, x'=0$ & $F_{0110}$ & $[M_0 \times \bar{Z}; Y \times Y]$ \\
$r=0, r'=0$ & $F_{0101}$ & $[M_0 \times M_0 ; Y \times Y]$ \\
\hline
\end{tabular}

These faces meet in the following corners. 

\begin{tabular}{|l|l|l|l|}
\hline
Arising from & corner & geometry & in faces \\
\hline
$x=0, r=0, x'=0$ & $C_{1110}$ & $Y \times \bar{Z}$ & $F_{1010}, F_{0110}$\\
$x=0, r=0, r'=0$ & $C_{1101}$ & $Y \times M_0$ & $F_{1001}, F_{0101}$\\
$x=0, x'=0, r'=0$ & $C_{1011}$ & $\bar{Z} \times Y$ & $F_{1010}, F_{1001}$\\
$r=0, x'=0, r'=0$ & $C_{0111}$ & $M_0 \times Y$ & $F_{0110}, F_{0101}$\\
$x=0, x'=0, r=0, r'=0$ blowup & $C_{1111}$ & $SN^+ (Y \times Y)$ & $F_{1010}, F_{1001}, F_{0110}, F_{0101}.$\\
\hline
\end{tabular}

To see that $\mathcal{D}$ is a smooth submanifold of $\mathcal{S}^2 _b$ we consider  the function 
\begin{equation} \label{eq:fp} f(p) = x(p) r(p) - x'(p) r'(p). \end{equation}
Away from the boundary faces $f$ is smooth with non-vanishing differential.  In a neighborhood of $S_{11} - F_{1001},$ let
$$f_{0110} (p) = \frac{x(p)}{x'(p)} - \frac{r'(p)}{r(p)}.$$
Since $x'r f_{0110} = f,$ we see that $f_{0110}$ is smooth near the $\epsilon =0$ boundary faces away from where those faces meet $F_{0110}.$  Moreover, wherever defined, $f_{0110}$ has nonvanishing differential and the zero set of $f_{0110}$ coincides with that of $f$ away from $F_{0110}.$  Similarly, let 
$$f_{1001} (p) = \frac{x(p)}{r'(p)} - \frac{x'(p)}{r(p)}.$$
$f_{1001}$ is smooth with nonvanishing differential and has the same vanishing set as $f$ in a neighborhood of $\{\epsilon =0 \} - F_{1001}.$  This shows that $\mathcal{D}$ is a smooth submanifold of $\mathcal{S}^2 _b.$  While the acc double space will not be used here, we note that the acc double space, with an additional blowup along the diagonal for $\epsilon \geq 0$, would be the natural space on which to study the resolvent behavior under ac convergence.  

\section{Analytic Preliminaries}

Since we are working on manifolds with singularities, corners and boundaries, we briefly review some key features of the analysis in these settings.  

\subsection{Polyhomogeneous conormal functions}

On manifolds $M$ with corners having a consistent local product structure near each boundary and corner, a natural class of functions (or sections) with good regularity near the boundary and corners are the polyhomogeneous conormal functions (or sections).  For a complete reference on polyhomogeneity on manifolds with corners, see \cite{edge}.  In a neighborhood of a corner, we have coordinates $(x_1, \ldots, x_k, y_1, \ldots, y_{n-k})$ where $x_1, \ldots, x_k$ vanish at this corner and $(y_1, \ldots , y_{n-k})$ are smooth local coordinates on a smooth compact $n-k$ manifold $Y.$  The edge tangent bundle $\mathcal{V}_e $ in a neighborhood of this corner is spanned over $\cC^{\infty}(M)$ by the vector fields, 
$$\{ x_i \partial_{x_i} ,  \partial_{y ^{\alpha}} \}.$$
The basic conormal space of sections is 
$$\mathcal{A}^0 (M_0 ) = \{ \phi :  V_1 ... V_l \phi \in L^{\infty} (M_0 ), \forall \, V_i \in \mathcal{V}_e , \textrm{ and } \forall \, l \}.$$
Let $\alpha$ and $p$ be multi indices with $\alpha_j \in \mathbb{C}$ and $p_j \in \mathbb{N}_0$.  Then we define
$$\mathcal{A}^{\alpha , p} (M_0 ) = x^{\alpha} (\log x)^p \mathcal{A}^0 .$$
The space $\mathcal{A}^* $ is the union of all these spaces, for all $\alpha$ and $p.$  The space $\mathcal{A}^* _{phg} (M_0 )$ consists of all conormal distributional sections which have an expansion of the form
$$\phi \sim \sum_{Re(\alpha_j ) \to \infty} \sum_{p=0} ^{p_j} x^{\alpha_j} (\log x)^p a_{j, p} (x, y), \, a_{j, p} \in \mathcal{C}^{\infty} .$$
We define an index set to be a discrete subset $E \subset \mathbb{C} \times \mathbb{N}_0$ such that 
$$(\alpha_j , p_j ) \in E, \quad |( \alpha_j , p_j ) | \to \infty \implies Re(\alpha_j ) \to \infty.$$
Then, the space $\mathcal{A}^E _{phg} (M_0)$ consists of those distributional sections $\phi \in \mathcal{A} ^* _{phg}$ having polyhomogeneous expansions with  $(\alpha_j , p_j ) \in E.$  

\subsection{Conic differential operators and $b$-operators}

Let $(M_0 , g_0 )$ be a Riemannian manifold with isolated conic singularity, defined by $x=0$, so that in a neighborhood of the singularity, 
$$M_0  \cong (0, x_1)_x \times Y, \quad g_0 = dx^2 + x^2 h(x) $$
with $h(x ) \to h$, where $(Y, h)$ is a smooth, $n-1$ dimensional compact manifold.  A conic differential operator of order $m$ is a smooth differential operator on $M_0$ such that in a neighborhood of the singularity it can be expressed 
$$A = x^{-m} \sum_{k=0} ^m B_k (x) (-x \partial_x )^k$$
with $B_k \in \mathcal{C}^{\infty} ( (0,x_1), \textrm{Diff}^{m-k} (Y))$, where Diff$^j (Y)$ denotes the space of differential operators of order $j \in \mathbb{N}_0$ on $Y$ with smooth coefficients.  The cone differential operators are elements of the cone operator calculus; for a detailed description, see \cite{lesch}.  These cone operators are closely related to $b$-operators.   A $b$-operator of order $m$ is a smooth differential operator such that near the boundary it can be expressed 
$$A =  \sum_{k=0} ^m B_k (x) (-x \partial_x )^k$$
with $B_k \in \mathcal{C}^{\infty} ( (0, x_1), \textrm{Diff}^{m-k} (Y))$.  We see that a cone differential operator of order $m$ is equal to a rescaled $b$-differential operator of order $m.$  In other words, if $A$ is an order $m$ cone differential operator then $x^m A$ is a $b$-differential operator.  In local coordinates $(x, y_1 , \ldots , y_{n-1} )$ near the boundary of $M_0,$ a $b$-operator may be expressed as 
$$A = \sum_{j + |\alpha| \leq m} a_{j, \alpha} (x, y) (-x \partial_x)^j (\partial_{y ^{\alpha} }).$$
The $b$-symbol of $A$ is 
$$^b \sigma_m (A) = \sum_{j + |\alpha| = m } a_{j, \alpha} (x, y) \lambda^{j} \eta^{\alpha} .$$
Here $\lambda$ and $\eta$ are linear functions on $^b T^* M_0  $ defined by the coordinates so that a generic element of $^b T^* M_0 $ is 
$$ \lambda \frac{dx}{x} + \sum_{i=1} ^n \eta_i dy_i .$$
The $b$-operator is $b$-elliptic if the symbol $^b \sigma_m (A) \neq 0$ on $^b T^* M_0 - \{0\}$.

The scalar Laplacian on $M_0$ is 
$$x^{-2} \{ (-x \partial_x )^2 + (-n + 1+ x H^{-1} (\partial_x H )(-x \partial_x )) + \Delta_{h(x)} \} = x^{-2} L_b$$
where $L_b$ is an elliptic order $2$ $b$-operator and $H$ is a smooth function depending on the metric.  Similarly, a geometric Laplacian $\Delta_0$ on $M_0$ is also of the form 
$$\Delta_0 = x^{-2} L_b,$$
for an elliptic order two $b$-operator acting on sections of the vector bundle.  The Schwartz kernel of $L_b$ is a distribution on the $b$-double space $M^2 _{0,b}.$  By the $b$-calculus theory, (see \cite{tapsit}) $L_b$ has a parametrix $G_b$, such that $G_b$ is a $b$-operator of order $-2$ with
$$G_b L_b  = I - R $$
where $I$ is the identity operator and $R$ is a $b$-operator with polyhomogeneous Schwartz kernel on the $b$-double space.  Then, for any $u \in \mathcal{L}^2 (x^{n-1} dx dy)$ with $\Delta_0 u = f \in \mathcal{L}^2 (x^{n-1} dx dy)$
$$ (x^{2} G_b ) (x^{-2}L_b  u ) = (x^2 G_b ) f = u - Ru \quad \implies u = x^2 G_b f + Ru = \alpha + \beta.$$
The first term, $\alpha \in x^2 H^2 _b \subset x^2 \mathcal{L}^2 (x^{n-1} dx dy)$.   The second term $\beta \in \mathcal{L}^2 (x^{n-1} dx dy)$ has a polyhomogeneous expansion as $x \to 0,$
$$ \beta \sim \sum_{j=0} ^{\infty} \sum_{k=0} ^{N_j} x^{\gamma_j + k} \varphi_j (y). $$ Above $\gamma_j$ is an indicial root for the operator $L_b$ and $\varphi_j$ is an eigensection for the induced geometric Laplacian on $(Y, h)$.  Then,  
\begin{equation} \label{eq:conphg} u = \alpha + \sum_{j=0} ^{\infty} \sum_{k=0} ^{N_j} x^{\gamma_j + k} \varphi_j (y) \end{equation}
where $\alpha \in x^2 \mathcal{L}^2 (x^{n-1} dx dy)$.   This decomposition plays a key role in the proof of spectral convergence.   

\subsection{Friedrich's domain of the conic Laplacian} 

A geometric Laplacian $\Delta_0$ on a  conic manifold is an unbounded operator on $\mathcal{L}^2$ sections of the bundle.  It can be extended to various domains in $\mathcal{L}^2;$ the minimal domain $\mathcal{D}_{\textrm{min}}$ is the $\cL^2$ closure of the graph of $\Delta_0$ over $\mathcal{C}^{\infty} _0 .$  The largest domain  $\mathcal{D}_{\textrm{max}}$ is the $\mathcal{L}^2$ closure of the graph of $\Delta_0$ over $\mathcal{L}^2.$  Each of these domains are dense in $\mathcal{L}^2 (M_0 ),$ and the extension of the Laplacian to either domain is a closed operator.  On complete manifolds $\mathcal{D}_{\textrm{min}} = \mathcal{D}_{\textrm{max}}$ by the Gaffney-Stokes Theorem \cite{ga}.  However, $M_0$ is incomplete and so for a  general geometric Laplacian these domains will not be equal.  The Friedrich's domain $\mathcal{D}_F$ lies between $\mathcal{D}_{\textrm{min}}$ and $\mathcal{D}_{\textrm{max}}$ and is the closure of the graph of $\Delta_0$ in $\mathcal{L}^2$ with respect to the densely defined Hermitian form,
$$ Q(u, v) = \int_{M_0} \langle \nabla u, \nabla v \rangle . $$
The extension of the Laplacian to the Friedrich's domain, known as the Friedrich's extension of the Laplacian, preserves the lower bound and is essentially self adjoint.  Here, we work exclusively with the Friedrich's extension of the Laplacian.        

For elements of $\mathcal{D}_{max}$, with $u \in \mathcal{L}^2$ and $\Delta_0 u = f \in \mathcal{L}^2$, we have the expansion (\ref{eq:conphg}) from the preceding section, 
$$ u = \alpha + \sum_{j=0} ^{\infty} \sum_{k=0} ^{N_j} x^{\gamma_j + k} \varphi_j (y). $$
The volume form on $M_0$ near the singularity is asymptotic to $x^{n-1} dx dy$.  Therefore, the exponents $\gamma_j$ must all be strictly greater than $-\frac{n}{2}$.  For $v \in \mathcal{D}_{\textrm{min}} \subset \mathcal{D}_{max}$ the decomposition (\ref{eq:conphg}) and the definition of $\mathcal{D}_{min}$ imply that $\mathcal{D}_{\textrm{min}} \subset  x^2 \mathcal{L}^2$.  The equality of $\mathcal{D}_{\textrm{min}}$ and $\mathcal{D}_{\textrm{max}}$ then depends on the indicial roots of $L_b = x^2 \Delta_0$.  For further discussion of domains of the conic Laplacian, see \cite{domains}, whose results include:
$$\mathcal{D}_F = \{ f \in \mathcal{L}^2 : \Delta_0 f \in \mathcal{L}^2 \textrm{ and } f = \mathcal{O}( x^{\frac{2-n + \delta }{2}}) \textrm{ as $x \to 0$, for some $\delta > 0$} \}.$$  
We will use this characterization of the domain of the (Friedrich's extension of the) Laplacian in the proof of the first theorem.  

\section{Spectral Convergence}

We now have all the necessary ingredients to prove spectral convergence.  

\begin{prethm2}
\label{prethm2}
Let $(M_0, g_0 )$ be a compact Riemannian $n$-manifold with isolated conic singularity,  and let $(Z, g_z)$ be an asymptotically conic space, with $n \geq 3.$  Assume $(M, g_{\epsilon} )$ converges asymptotically conically to $(M_0 , g_0).$  Let $(E_0, \nabla_0)$ and $(E_z, \nabla_z)$ be Hermitian vector bundles over $(M_0, g_0)$ and $(Z, g_z),$ respectively, so that each of these bundles in a neighborhood of the boundary is the pullback from a bundle over the cross section $(Y, h).$  Let $\Delta_0,$ $\Delta_z$ be the corresponding Friedrich's extensions of geometric Laplacians, and let $\Delta_{\epsilon}$ be the induced geometric Laplacian on $(M, g_{\epsilon}).$  Assume $\Delta_z$ has no $\mathcal{L}^{2}$ nullspace.  Then the accumulation points of the spectrum of $\Delta_{\epsilon}$ as $\epsilon \to 0$ are precisely the points of the spectrum of $\Delta_0,$ counting multiplicity.   
\end{prethm2}

The theorem follows from the following three statements:  the inclusion $accumulation$  $\sigma ( \Delta_{\epsilon} )  \subset \sigma ( \Delta_{0}), $ the reverse inclusion $accumulation$ $ \sigma ( \Delta_{\epsilon} )  \supset \sigma ( \Delta_{0}), $ and correct multiplicities.

\subsection{Accumulation $\sigma ( \Delta_{\epsilon} ) \subset \sigma ( \Delta_{0}) $}

We extract a smoothly convergent sequence of eigensections corresponding to a converging sequence of eigenvalues as $\epsilon \to 0,$ and show that the limit section of this sequence is an eigensection for the conic metric and its eigenvalue is the accumulation point.  For this argument, we work with sequences of metrics $\{ g_{\epsilon_j} \}$ which we abbreviate $\{ g_j  \}$ with Laplacians $\Delta_j.$     

Let $\lambda(\epsilon_j )$ be an eigenvalue of $\Delta_j$, with eigensection
$f_j$.  Assume that $\lambda(\epsilon_j) \to \bar{\lambda}$. Over any compact
set $K \subset M_0 ^0$, the metric $g_j = g_{\epsilon_j}$ converges smoothly to $g_0$ by Lemma 1,
thus so do the coefficients of $\Delta_j - \Delta(\epsilon_j)$. Hence,
normalizing $f_j$ by $\sup_{M} |f_j| =1$, it follows using standard elliptic
estimates and the Arzela-Ascoli theorem that $f_j$ converges in ${\mathcal C}^\infty$
on any compact subset of $M_0 ^0$. Furthermore the limit section $\bar{f}$ satisfies the limiting equation
$$ \Delta_0 \bar{f} = \bar{\lambda} \bar{f}.  $$
However, we do not know yet that $\bar{f} \not\equiv 0$, nor, even if this limit is nontrivial, that it lies in the domain of the Friedrichs extension of $\Delta_0.$  This is the content of the arguments below.  

\subsubsection{Weight Functions} 
Let $\phi_{\epsilon} : M_{0, \epsilon} - M_{0, \delta} \to Z_{1/\epsilon} - Z_{1/\delta}$ as in definition \ref{resblow}.  We identify $Z_{1/ \delta}$ with a fixed $K \subset U \subset M$ so that $M_{0, \epsilon} - M_{0, \delta} \cong (U-K),$ $K \cong Z_{1/\delta}.$  

Let
\begin{displaymath}
w_{\epsilon} = \left \{ \begin{array}{ll} c & \textrm{ on $M-U$.}\\ 
\epsilon \, (\phi_{\epsilon}^{-1})^* \rho & \textrm{ on $U-K$.} \\
c\epsilon & \textrm{ on $K$.} \end{array} \right. \end{displaymath}

Above, $c$ is a constant and no generality is lost by assuming $c=1.$  Let $w_j = w_{\epsilon_j}$.  For some $\delta >0$ to be chosen later, replacing $f_j$ by $\frac{f_j}{|| w_j ^{\delta} f_j ||_{\infty}}$ we assume the supremum of $|f_j w_j ^{\delta} |$ is $1$ on $M$.  Since $M$ is compact, $|f_j|$ attains a maximum at some point $p_j \in M$, and we may assume $p_j$ converges to some $\bar{p} \in M$.  The argument splits into three cases depending on how and where $p_j$ accumulates in $M.$    

\subsubsection{Case 1: $w_j (p_j ) \to c >0 $ as $j \to \infty$.} 

In this case, the points $\{p_j \}$ are accumulating in a compact subset of $M_0 ^0$ at some point $\bar{p} \neq p$ (the singularity).  The maximum of $| f_j w_j ^{\delta} |$ on $M$ is 1 and occurs at $p_j$ so   
$$|f_j | \leq w_j ^{-\delta} \textrm{ on $M$ for each $j$ } \implies |f_j (p_j) | \to c^{-\delta} \textrm{ as } j \to \infty.$$

The locally uniform $\mathcal{C}^{\infty}$ convergence of $f_j$ to $\bar{f}$ implies that $|\bar{f}|$ satisfies a similar bound,
$$|\bar{f} | \leq x^{-\delta} \textrm{ as $x \to 0,$}$$
and clearly $|\bar{f} (\bar{p})| = c^{- \delta} \neq 0.$  
By the dimension assumption $n \geq 3$ and the characterization of the Friedrich's Domain of the Laplacian, we may choose $\delta$ so that 
$$\frac{2-n}{2} < -\delta < 0.$$
Then $\bar{f}$ lies in the Friedrich's domain of the Laplacian and satisfies
$$\Delta_0 \bar{f} = \bar{\lambda} \bar{f},$$
so $\bar{\lambda}$ is an eigenvalue of $\Delta_0.$

\subsubsection{Case 2:  $| w_j (p_j ) | \leq c (\epsilon_j  )$ as $j \to \infty$.} 

Analysis on $Z$ in this case leads to a contradiction.  Let $\phi_j = \phi_{\epsilon_j},$ $\tilde{f_j} = f_j(\phi_j ^{-1}).$   Let $\tilde{p_j} = \phi_j (p_j)$.  Because $| f_j w_j ^{\delta} |$ attains its maximum value of $1$ at $p_j$, $|\tilde{f_j} (\tilde{p_j} )| = (w_j (p_j) )^{-\delta}$.  Rescale $f_j$ and $\tilde{f_j}$, replacing them respectively with $(w_j (p_j) )^{\delta} f_j$ and $(w_j (p_j) ) ^{\delta} \tilde{f_j} $ so that the maximum of $| \tilde{f_j} \rho^{\delta} |$ occurs at the point $\tilde{p_j} \in Z_j$ and is equal to $1$.  Since $w_j (p_j) = \mathcal{O} (\epsilon_j  )$, $\rho (\tilde{p_j} ) = \epsilon_j ^{-1} w_j (p_j ) $ stays bounded for all $j$, and so we assume $\tilde{p_j}$ converges to $\tilde{p} \in Z$.  By Lemma 1, $(Z_j , \epsilon_j ^{-2} \phi_j ^* g_j |_U )$ converges smoothly to $(Z_j , g_Z  )$.  This implies the following equation is satisfied by $\tilde{f_j}$ on $Z_j,$
$$\Delta_Z \tilde{f_j} = \epsilon_j ^2 \lambda(\epsilon_j ) \tilde{f_j} + O(\epsilon_j).$$
Since the $\lambda(\epsilon_j )$ are converging to $\bar{\lambda}$  and $| \tilde{f_j} \rho^{\delta} | \leq 1$ on $Z_j$, we have 
$$\Delta_Z \tilde{f_j} \to 0 \textrm{ as $j \to \infty,$ on any compact subset of $Z$.}$$
This implies $f_j \to \bar{f}$ on $M$ and correspondingly, $\tilde{f_j} \to \tilde{f}$ locally uniformly $\mathcal{C}^{\infty}$ on $Z$ and $\tilde{f}$ satisfies
$$\Delta_Z \tilde{f} = 0, \quad |\tilde{f}  \rho^{\delta}| \leq 1,$$
where equality holds in the second equation at the point $p.$  This shows that $\tilde{f}$ is not identically zero on $Z$ and $\tilde{f} = O(\rho^{-\delta})$ as $\rho \to \infty.$  Since $\tilde{f}$ is smooth on any compact subset of $Z$ and is therefore in $\cL^2 _{loc} (Z),$ choosing $\delta > n-2$ contradicts the assumption that $Z$ has no $\cL^2$ nullspace.  

\subsubsection{Case 3:  $ w_j (p_j) \to 0, \quad \frac{\epsilon_j}{w_j (p_j )} \to 0$ as $j \to \infty$.} 

In this case, the points $\phi_j (p_j) \to \infty$ in $Z,$ so we rescale and derive a contradiction on the complete cone over $(Y, h)$.  Consider the coordinates $(\rho, y)$ on $Z$ defined for $\rho \geq \rho_1$.  In these coordinates $g_z = d\rho^2 + \rho^2 h(\rho)$.  Let $r_j = \frac{\epsilon_j }{w_j (p_j) } \rho$ and $\tilde{g_j} $ on $Z_j$ be defined by 
$$ \tilde{g_j} = \left( \frac{ \epsilon _j }{ w_j (p_j) } \right)^2 g_Z . $$
Then, 
$$(Z_j, \tilde{g_j} ) \cong \left( \left(\frac{ \rho_1 \epsilon_j}{w_j (p_j)}, \frac{1}{w_j (p_j)} \right) \times Y, \, \, dr_j ^2 + r_j ^2 h( r_j \epsilon_j / w_j (p_j) \right). $$
As $j \to \infty,$ $h(r_j \epsilon_j / w_j (p_j))$ converges smoothly to $h$, and 
$$\tilde{g_j} \to g_C = dr^2 + r^2 h$$
on the complete cone $C$ over $(Y, h).$  Let $\tilde{f_j} = w_j (p_j)^{\delta} f_j(\phi_j^{-1}).$   Since $|f_j w_j ^{\delta} | \leq 1$ with equality at $p_j ,$ $$|\tilde{f_j} r_j ^{\delta} | \leq 1 \textrm{ on $(Z_j , \tilde{g_j} )$ with equality at } \tilde{p_j} = \phi_j (p_j).$$
Let $\tilde{\Delta_j}$ on $Z_j$ be the Laplacian induced by $\tilde{g_j}$ on $Z_j,$
$$\tilde{\Delta_j} = \frac{ w_j (p_j)^2}{\epsilon^2} \Delta_Z,$$ 
so
$$\tilde{\Delta_j} \tilde{f_j} = \epsilon_j ^2 \frac{w_j(p_j) ^2}{\epsilon_j^2}   \lambda (\epsilon_j ) \tilde{f_j} + O (\epsilon_j) \textrm{ on $(Z_j , \tilde{g_j} )$.}$$
Since $w_j (p_j) \to 0$ as $j \to \infty,$ there is a locally uniform $\mathcal{C}^{\infty}$ limit $f_c$ of $\{\tilde{f_j} \}$ on $C$ which satisfies 
$$|f_c r^{\delta} | \leq 1, \quad \Delta_c f_c = 0.$$
Since the points $\tilde{p_j}$ stay at a bounded radial distance with respect to the radial variable $r_j$ on $Z_j$, we may assume $\tilde{p_j} \to p_c$ for some $p_c \in C$.  At this point, $|f_c (p_c) r(p_c )^{\delta} | = 1$ so $f_c$ is not identically zero.  By separation of variables (see, for example, \cite{lo}), $f_c$ has an expansion in an orthonormal eigenbasis $\{ \phi_j \}$ of $\mathcal{L}^2 (Y, h),$
$$f_c = \sum_{j \geq 0} a_{j, + } r^{\gamma_{j, +}} \phi_j (y) + a_{j, -} r^{\gamma_{j, -}} \phi_j (y) $$
where $\gamma_{j, +/-}$ are indicial roots corresponding to $\phi_j$ and $a_{j, +/-} \in \mathbb{C}$.  In order for $|f_c r^{\delta} | \leq 1$ globally on $C$, we must have only one term in this expansion, $f_c = a_j r^{- \delta} \phi_j (y)$.  Because the indicial roots are discrete, we may choose $\delta$ so that $-\delta$ is not an indicial root.  This is a contradiction.

\subsection{$\sigma (\Delta_0 ) \subset $ Accumulation $\sigma(\Delta_{\epsilon}) $ } 

We use the Rayleigh-Ritz characterization of the eigenvalues \cite{chavel}.  Let $\lambda_l (\epsilon_j )$ be the $l^{th}$ eigenvalue of $\Delta_j$ and let 
$$ R_j (f) := \frac{\langle \triangledown f, \triangledown f \rangle_{j} }{\langle f, f \rangle_{j}}. $$

The subscript $j$ indicates that the inner product is taken with respect to the $\mathcal{L}^2$ norm on $M$ with the $g_j$ metric.  The eigenvalues are characterized using Mini-Max by
$$ \lambda_l (\epsilon_j ) = \textrm{inf}_{\textrm{dim $L=l$, } L \subset \mathcal{C}^1 (M)} \textrm{  sup}_{f \in L, f \neq 0} \, R_j(f).$$
Similarly this characterization holds for the eigenvalues of the (Friedrich's extension of the) conic Laplacian which are known to be discrete (see \cite{chee}, for example).  Because $\mathcal{C}^{\infty} _0 (M_0 )$ is dense in $\mathcal{L}^2 (M_0 )$ we may restrict to subspaces contained in $\mathcal{C}^{\infty} _0 (M_0 )$.  Then, the $l^{th}$ eigenvalue of $\Delta_0$ is
$$ \bar{\lambda}_l = \textrm{inf}_{\textrm{dim $L=l$, } L \subset \,  \mathcal{C}^{\infty} _0 (M_0 )} \textrm{  sup}_{f \in L, f \neq 0} \, R_0 (f). $$

Let $\bar{\lambda}_l$ be the $l^{th}$ eigenvalue in the spectrum of $\Delta_0$.  Fix $\epsilon >0$.  Then there exists $L \subset \mathcal{C}^{\infty} _0$ with $dim(L) = l$ and  
$$ \textrm{  sup}_{f \in L, f \neq 0} \, R_0 (f) < \bar{\lambda}_l + \epsilon.$$
Since any $f \in L$ is also in $\mathcal{C}^{\infty} _0 (M)$ and because $L$ is finite dimensional, by the local convergence of $g_j$ to $g_0,$ for large $j$
$$|R_j (f) - R_0 (f) | < \epsilon \textrm{ for any $f \in L$.}$$
Since $\lambda_l (\epsilon_j )$ is the infimum 
$$\lambda_l (\epsilon_j ) \leq \bar{\lambda}_l + 2 \epsilon.$$
This shows $\{ \lambda_l (\epsilon_j ) \}$ is bounded in $j,$ and so we extract a convergent subsequence and a corresponding convergent sequence of eigensections which exists by the preceding arguments.  For each $l$ we take 
$$\lambda_l (\epsilon_j ) \to \mu_l \leq \bar{\lambda}_l,$$
$$f_{j, l} \to u_l , \quad \Delta_0 u_l = \mu_l u_l .$$
These limit eigensections $u_l$ are seen to be orthogonal as follows.  Fix $l,$ $k,$ with $f_{j, k} \to u_k$ and $f_{j, l} \to u_l$.  Since $\mathcal{C}^{\infty} _0 (M_0 )$ is dense in $\mathcal{L}^2 (M_0 )$ we may choose a smooth cutoff function $\chi$ vanishing identically near the singularity in $M_0$ such that 
$$|| \chi u_k - u_k ||_{L^2 (M_0 )} < \epsilon,$$
$$|| \chi u_l - u_l ||_{L^2 (M_0 )} < \epsilon,$$
$$Vol_j (M \, - \, supp(\chi)) < \epsilon.$$
Then on the support of $\chi$, $g_j \to g_0$ uniformly so for large $j,$
$$|\langle u_k , u_l \rangle_0 - \langle \chi u_k , \chi u_l \rangle_0 | < \epsilon,$$
$$|\langle \chi u_k , \chi u_l \rangle_0 - \langle \chi u_k , \chi u_l \rangle_j | < \epsilon,$$
$$|\langle \chi u_k , \chi u_l \rangle_j - \langle \chi u_k , \chi f_{j, l} \rangle_j | < \epsilon,$$
$$|\langle \chi u_k , \chi f_{j,  l} \rangle_j - \langle \chi f_{j, k} , \chi f_{j, l} \rangle_j | < \epsilon.$$
Since the eigensections for $\Delta_j$ were chosen to be orthonormal and the volume of $\left( M - \textrm{support}(\chi) \right)$ is small with respect to $g_j$,
$$|\langle \chi f_{j, k} , \chi f_{j, l} \rangle_j | < 2 \epsilon.$$
Thus, $\langle u_k , u_l \rangle_0$ can be made arbitrarily small and $u_k, u_l$ are orthogonal for $l \neq k.$  We complete this basis to form an eigenbasis of $\mathcal{L}^2 (M_0 )$.  Let $\bar{f_l}$ be an arbitrary element of this eigenbasis, with eigenvalue $\bar{\lambda_l}$.  We wish to show that this $\bar{f_l}$ is actually the $u_l$ above, defined to be the limit of (a subsequence of) $\{f_{j, l} \}$, and hence the corresponding $\mu_l$ is equal to $\bar{\lambda_l}$.  Again, assume the smooth cut-off function $\chi$ is chosen so that 
$$|| \chi \bar{f_l} - \bar{f_l} ||_{L^2 (M_0 ) } < \epsilon .$$
For each $j$ we expand $\chi \bar{f_l}$ in eigensections of $\Delta_j,$
$$\chi \bar{f_l} = \sum_{k = 0} ^{\infty} a_{j, k} f_{j, k},  \quad \textrm{where }a_{j, k} = \langle \chi \bar{f_l} , f_{j, k} \rangle_j .$$
Now, fix $k$ and choose $\chi$ such that 
$$|| \chi u_k - u_k ||_{L^2 (M_0 )} < \epsilon. $$
Then, 
$$|\langle \chi \bar{f_l}, f_{j, k} \rangle_0 - \langle \chi \bar{f_l} , f_{j, k} \rangle_j | < \epsilon,$$
$$|\langle \chi \bar{f_l} , f_{j, k} \rangle_j - \langle \chi \bar{f_l} , u_k \rangle_j | < \epsilon,$$
$$| \langle \langle \chi \bar{f_l} , u_k \rangle_j - \langle \chi \bar{f_l} , u_k \rangle_0 | < \epsilon,$$
$$| \langle \chi \bar{f_l} , u_k \rangle_0 - \langle \bar{f_l} , u_k \rangle_0 | < \epsilon.$$
By the orthogonality $\langle \bar{f_l} , u_k \rangle_0 = 0$ if $\bar{f_l} \neq u_k$, and otherwise is $1,$ so for each $k$, $a_{j, k} \to 0$ as $j \to \infty,$ for all $k$ with $u_k \neq \bar{f_l}$.  Because $f_l$ is not identically zero there must be some $k$ with $u_k = \bar{f_l}$.  This shows that \em every \em eigensection of $\Delta_0$ is the limit of (a subsequence of) $\{f_{j, k} \}$ and the corresponding eigenvalue $\bar{\lambda_k}$ is the limit of the corresponding eigenvalues. 

\subsection{Correct Multiplicities}

We argue here by contradiction.  Let $\lambda$ be an eigenvalue for $\Delta_0$ with $k$ dimensional eigenspace spanned by $u_1 , \ldots, u_k$.  Assume $\lambda$ occurs as an accumulation point of multiplicity less than $k,$ without loss of generality assume multiplicity $k-1.$  However, preceding arguments imply the existence of a subsequence $\{f_{k,j} \}$ of $\Delta_j$ with $f_{k,j} \to u_k,$ which shows that $\lambda$ is achieved as an accumulation point of multiplicity at least $k.$  Conversely, assume $\lambda$ occurs as an accumulation point of multiplicity $k+1.$  By preceding orthogonality argument, the limit section $u_{k+1}$ of the converging sequence $f_{k+1, j}$ is orthogonal to $\{u_1, \ldots, u_k \}.$  This is a contradiction.  
\heart

\section{Heat Kernels}

The heat kernels for each of the geometries in ac convergence are elements of a pseudodifferential heat operator calculus that is defined on the corresponding heat space.  For the details in the construction of these heat calculi, kernels and spaces, see \cite{moi}.

\subsection{$b$-heat kernel}

Let $(M, g)$ be a $b$-manifold with local coordinates $z=(x, y)$ in a neighborhood of $\partial M$ so that near the boundary 
$$g = \frac{dx^2}{x^2} + h(x, y).$$
Let $(z, z')$ be coordinates on $M \times M$ and let $\Delta_b$ be a geometric Laplacian on $M.$  The $b$-heat kernel $H(z, z', t)$ is the Schwartz kernel of the fundamental solution of the heat operator $\partial_t + \Delta_b.$  The heat kernel is a distributional section that acts on smooth sections of $M,$ and satisfies
$$(\partial_t + \Delta_b )H(z, z', t) = 0, \, t>0,$$
$$H|_{t=0} = \delta(z-z'),$$
and by self adjointness since we work with the Friedrich's extension of $\Delta_b$, 
$$H(z, z', t) = H(z', z, t)^*.$$
For a smooth section $u$ on $M,$ 
$$u(z,t) := \int_M \langle u(z'), H(z, z', t) \rangle dz'$$
satisfies
$$( \partial_t + \Delta_b ) u(z, t) = 0 \textrm{ for } t>0, \quad u(z, 0) = u(z).$$
Physically, $u(z, t)$ describes the heat on $M$ at time $t>0$ where the initial heat applied to $M$ is given by $u(z).$  

Recall the Euclidean heat kernel, 
$$G(z, z', t) = (4 \pi t)^{-\frac{n}{2}} \exp\left( -\frac{|z-z'|^2}{2t} \right).$$  
For a compact manifold without boundary, the heat kernel can be constructed locally using the Euclidean heat kernel and geodesic normal coordinates \cite{rose}.  On the interior of a manifold with boundary (or singularity) the Euclidean heat kernel is also a good model, however, near the boundary (singularity) a different construction is required.

\subsubsection{$b$-heat space}

It is convenient to study the heat kernel on a manifold with boundary (or singularity) as an element of a heat operator calculus defined on a corresponding heat space.  This space is a manifold with corners constructed from $M \times M \times \R^+$ by blowing up along submanifolds at which the heat kernel may have interesting or singular behavior.  For example, the diagonal is always blown up at $t=0,$ since away from the boundary the heat kernel behaves like the Euclidean heat kernel which is singular along the diagonal at $t=0.$  For the $b$-heat space we first blow up the codimension 2 corner at the boundary in both copies of $M.$  The \em $b$-heat space \em $M^2 _{b,h}$ is then,  
$$M^2 _{b, h} = \left[ M^2 _b \times \mathbb{R}^+ _t ; \Delta(M \times M) \times \{t=0\}, dt \right],$$
where $\Delta(M \times M)$ is the diagonal in $M \times M$ and $M^2 _b$ is the $b$-double space (\ref{bdouble}).   The $b$-heat space has five boundary faces, two of which result from blowing up.  The remaining three boundary faces are at $t=0$ off the diagonal and at the boundary in each copy of $M.$  More precisely, we have the following.\footnote{The subscript ``$d$'' indicates a face created by blowing up along the diagonal, so for example $F_{d2}$ is the face created by blowing up along the diagonal where the scalar variable $t$ vanishes to second order.}  

\begin{tabular}{|l|l|l|}
\hline
Face & Geometry of face & Defining function in  local coordinates \\
\hline
$F_{110}$ & $N^+ ( Y \times Y) \times \R^+$ & $\rho_{110} = (x^2 + (x')^2)^{\frac{1}{2}}$ \\
$F_{d2}$ & $PN^+_t ( \Delta(M \times M))$ & $\rho_{d2} = (|z-z'|^4 + t^2)^{\frac{1}{4}}$ \\
$F_{100}$ & $Y \times (M - \partial M) \times \R^+$ & $\rho_{100} = x$\\
$F_{010}$ & $Y \times(M - \partial M) \times  \R^+$ & $\rho_{010} = x'$\\
$F_{001}$ & $(M -  \partial M)^2 - \Delta(M\times M)$ & $\rho_{001} = t$\\
\hline
\end{tabular}

Above $PN^+ _t$ denotes the inward pointing $t$ parabolic normal bundle while $SN^+$ denotes the inward pointing spherical normal bundle.  Note that the local coordinates $x, x', t$ lift from $M \times M \times \R^+$ to $M^2 _{b,h}$ as follows,
$$\beta^* (x) = \rho_{110} \rho_{100}, \, \beta^* (x') = \rho_{110} \rho_{010}, \, \beta^* (t) = \rho_{d2}^2 \rho_{001},$$
so these coordinates are only \em local \em defining functions.  

\subsubsection{$b$-heat calculus}

The $b$-heat calculus consists of distributional section half density kernels on $M^2 _+ = M \times M \times \R^+$ which are smooth on the interior and lift to be polyhomogeneous on $M^2 _{b,h}$ with specified leading orders at the boundary faces.  By constructing the $b$-heat kernel as an element of the $b$-heat calculus, it is polyhomogenous on $M^2 _{b,h.}$ and we know the leading order terms.  Consequently, the $b$-heat kernel is polyhomogeneous at the boundaries and corners of $M^2 _+$ with specified leading orders.    Once the calculus is defined and the composition rule is proven, construction of the heat kernel as an element of the heat calculus is similar to solving an ordinary differential equation using Taylor series.  The following definition is from \cite{tapsit}.  

\begin{defbheat}
For any $k \in \mathbb{R}$ and index set $E_{110},$ $A$ is an element of  the $b$-heat calculus, $\Psi^{E_{110},  k} _{b,H} $ if the following hold.  
\begin{enumerate}
\item $A \in \mathcal{A}^{-\frac{1}{2} + E_{110}} _{phg} (F_{110} ).$ 
\item $A$ vanishes to infinite order at $F_{001},$ $F_{100},$ and $F_{010}.$
\item $A \in \rho_{d2} ^{- \frac{n+3}{2} -k} \mathcal{C}^{\infty} (F_{d2} ).$
\end{enumerate}
\end{defbheat}

Because the heat calculus is defined with half densities, the normalizing factors at $F_{110}$ and $F_{d2}$ simplify the composition rule.  An element $A$ of the $b$-heat calculus is the Schwartz kernel of an operator acting on a smooth half density section $f$  of $M$ by
$$A f (z, t) = \int_M \langle A (z, z', t), f(z') \rangle dz'.$$
Furthermore, $A$ acts by convolution in the $t$ variable so for a smooth half density section $f$ of $M \times \R^+ _t,$ 
$$A f(z, t) = \int_0 ^t \int_M \langle A(z, z', t-s), f(z', s) \rangle dz' ds.$$ 
Two elements of the $b$-heat calculus compose as follows.

\begin{bheatcomp}
Let $A \in \Psi^{k_a, \mathcal{A}} _{b, H}$ and let $B \in \Psi^{k_b, \mathcal{B}} _{b, H}.$  Then the composition, $A \circ B$ is an element of $\Psi^{k_a + k_b ,  \mathcal{A} + \mathcal{B}} _{b, H}.$
\end{bheatcomp}
The proof of this composition rule is in \cite{tapsit}.

\subsubsection{Construction of the $b$-heat kernel}

First we construct a model heat kernel $H_1$ as an element of the $b$-heat calculus that solves the heat equation up to an error vanishing to positive order at the boundary faces of $M^2 _{b,h}.$   On the interior of $M^2 _{b,h}$ restricting to a coordinate patch with coordinates $(z, z', t),$ we locally define 
$$H_1 (z, z', t) : = (4 \pi t)^{-n/2} e^{ ( |z-z'|_g )^2 /2t},$$
where $|  z-z' |_g$ is the distance from $z$ to $z'$ with respect to the metric $g.$  As $t \to 0$ away from the diagonal this construction immediately implies infinite order vanishing at $F_{001}.$  At $F_{d2}$ we solve exactly:  for each $p \in M$ and for each point $z \in F_{d2}$ in the fiber over $(p, p, 0)$ the heat kernel at that point is determined by the coefficients of the metric (and its derivatives) at $p.$  The \em normal operator \em of $\partial_t + \Delta_b $ is the restriction to $F_{110}$ of the lift of $\partial_t + \Delta_b$ to $M^2 _{b,h}.$   $H_1$ is defined at $F_{110}$ to be the kernel of a first order parametrix of this normal operator and is smooth at this face.  At $F_{100},$ $F_{010},$ and $F_{001}$ the model kernel vanishes to infinite order.  As constructed, $H_1$ satisfies
$$(\partial_t + \Delta) H_1 = K_1, \quad H_1 \in \Psi^{-2, 0} _{b, H}$$
where $K_1$ now vanishes to positive order at the boundary faces of $M^2 _{b, h}.$  Then, define
$$H_2 = H_1 - H_1 * K_1,$$
where now the error term 
$$K_2 = (\partial_t + \Delta)H_2$$
vanishes to higher order on the boundary faces of $M^2 _{b, h}$ by the composition rule.  This construction is iterated and Borel summation (see \cite{shu}) gives 
$H_{\infty} \in \Psi^{-2, 0} _{b, H}$ with $H_{\infty} - H_N = O(t^{N - \frac{n+3}{2}} ),$ for $N>0,$ so that 
$$(\partial_t + \Delta) H_{\infty} = K,$$
where $K$ vanishes to infinite order on the boundary faces of $M^{2} _{b, h}$ so we may push $K$ forward to $M \times M \times \R^+.$  We solve away the residual error term using the action of elements of the $b$-heat calculus as $t$-convolution operators.  As a $t$-convolution operator, the heat kernel is the identity.  Above, $K$ as a $t$-convolution operator is of the form $K = Id - A$ where $A$ is a Volterra operator and $Id$ is the identity.  An operator of this form has an inverse of the same form so defining 
$$H := H_{\infty} (Id - A)^{-1}$$
solves away this residual error term.  By construction the leading order behavior of the $b$-heat kernel is that of the model heat kernel and is summarized below.  

\begin{tabular}{|l|l|}
\hline
Face & Leading order\\ 
\hline
$F_{110}$ & $0;$ $O(t^{-\frac{1}{2}})$ as $t \to \infty.$ \\
$F_{d2}$ & $-\frac{n+3}{2} - (-2)$ \\
$F_{100}$ & $ \infty$ \\
$F_{010}$ & $\infty$ \\
$F_{001}$ & $\infty$ \\
\hline 
\end{tabular}

\subsection{Conic heat kernel}

Let $(M_0 , g_0 )$ be a compact manifold with isolated conic singularity and let $(E_0 , \nabla_0 )$ be a Hermitian vector bundle over $(M_0 , g_0 ).$  Let $\Delta_0$ be the Friedrich's extension of a geometric Laplacian on $(M_0 , g_0 )$ and let $Y$ be the smooth $n-1$ dimensional cross section of $M_0$ so that $\partial M_0 =Y.$  The conic heat kernel is constructed analogously to the $b$-heat kernel.  

\subsubsection{The conic heat space}

This construction comes from \cite{moo}.  The conic heat space $M^2 _{0, h}$ is a manifold with corners obtained from $M_0 \times M_0 \times \mathbb{R}^+ = M^2 _{0,+}$ by blowing up along two submanifolds,
$$M^2 _{0, h} := \left[ [M_0 \times M_0 \times \R^+ ; \partial M_0 \times \partial M_0 \times \{t=0\} , dt ] ; \Delta(M^0 _0 \times M^0 _0) \times \{t=0\}, dt \right].$$

The conic heat space has five boundary faces described in the following table in which $z= (x, y)$ and $z' = (x', y')$ are local coordinates in a neighborhood of the singularity in each copy of $M_0$ so that $x=0, x'=0$ define the singularity as well as the boundary of $M_0.$

\begin{tabular}{|l|l|l|}
\hline
Face & Geometry of face &  Defining function in local coordinates \\
\hline
$F_{112}$ & $PN^+ _t (Y \times Y)$ & $\rho_{112} = (x^4 + (x')^4 + t^2)^{\frac{1}{4}}$ \\
$F_{d2}$ & $PN^+ _t (\Delta(M^0 _0 \times M_0 ^0))$ & $\rho_{d2} = ( |z-z'|^4 + t^2)^{\frac{1}{4}}$ \\
$F_{100}$ & $Y \times \R^+ $ & $\rho_{100} = x$\\
$F_{010}$ & $ \R^+ \times Y$ & $\rho_{010} = x'$\\
$F_{001}$ & $M_0^0 \times M^0 _0 - \Delta(M^0 _0 \times M^0 _0)$ & $\rho_{001} = t$\\
\hline
\end{tabular}

Note that the coordinates $x, x', t$ lift from $M_0 \times M_0 \times \R^+$ to $M^2 _{0, h}$ as follows, 
$$\beta^* (x) = \rho_{100} \rho_{112}, \, \beta^* (x') = \rho_{010} \rho_{112}, \, \beta^* (t) = \rho_{112} ^2 \rho_{d2} ^2 \rho_{001},$$
so again these are only \em local \em defining functions.  

\subsubsection{The conic heat calculus}
Let $\mu$ be a conic half density on $M^2 _{0, +};$  we may assume 
$$\mu =(xx')^{\frac{n-1}{2}} \sqrt{\textrm{d}z \textrm{d}z' \textrm{d}t} = \sqrt{\textrm{d}V_c \textrm{d}t}.$$  
Fix also a smooth, nonvanishing half density, $\nu,$ on $M_{0,h} ^2.$  Elements of the conic heat calculus are distributional section half densities on $M^2 _{0,+}$ which are smooth on the interior and lift to be polyhomogeneous on $M^2 _{0,h}.$  

\begin{heatcalc}
Let $k \in \mathbb{R}$ and $E_{100}$ $E_{010}$ $E_{112}$ be index sets.  Then $A \in \Psi^{k, E_{100}, E_{010}, E_{112}} _{0, H}$ if the following hold.   
\begin{enumerate}
\item $A \in \mathcal{A}^{E_{100}} _{phg}$ at $F_{100}.$ 
\item $A \in \mathcal{A}^{E_{010}} _{phg}$ at $F_{010} .$ 
\item $A \in \mathcal{A}^{E_{112}} _{phg}$ at $F_{112} .$
\item $A$ vanishes to infinite order at $F_{001}.$
\item $A \in \rho_{d2} ^{- \frac{n+3}{2} -k} \mathcal{C}^{\infty} (F_{d2} ).$
\end{enumerate}
\end{heatcalc}

With this normalization the conic heat kernel has order $k =-2$ and the composition rule is the following.

\begin{composition}
Let $A \in \Psi^{A_{100} , A_{010} , A_{112} , k_a} _{0,h}$, and $B \in \Psi^{B_{100}, B_{010} , B_{112} , k_b} _{0, H}$ with the leading index terms satisfying
$$\beta_{112} + \alpha_{010}  > 0, \, \,  \alpha_{112} + \beta_{100}  > 0, \, \, -k_a >0, \, \,  -k_b >0, \, \,  \beta_{100} + \alpha_{010} > -1. $$

Then, the composition $B \circ A$ is an element of $\Psi^{A_{100}, B_{010}, \Gamma_{112} , k} _{0, H}$ with $\Gamma_{112} = A_{112} + B_{112}$
and $k = (k_a + k_b)$.
\end{composition}  

The proof of this theorem is in \cite{moi}, and is originally due to \cite{moo}, see also \cite{lo}, \cite{gil}.

The conic heat kernel is constructed analogously to the $b$-heat kernel first using a model heat kernel and then using the composition rule to iteratively solve away the error term.  Away from $F_{112}, F_{100}$ and $F_{010}$ the model heat kernel comes from the standard local construction using the Euclidean heat kernel.  At $F_{112}$ the model heat kernel comes from an explicit construction of the heat kernel for an exact cone transplanted to $F_{112}$ and extended smoothly to $F_{100}$ and $F_{010}.$  See, for example \cite{chee0} and \cite{moo} or for the case of the scalar Laplacian, \cite{moi}.  The boundary behavior of the conic heat kernel on $M^2 _{0, h}$ is that of the model heat kernel and is summarized in the following table.  

\begin{tabular}{|l|l|}
\hline
Face & Leading Order \\
\hline
$F_{112}$ & $-\frac{3}{2} + \frac{2n+1}{2} + 2 \mu_0$\\
$F_{100}$ & $\frac{-n-1}{2} + \mu_0$\\
$F_{010}$ & $\frac{-n-1}{2} + \mu_0$\\
$F_{d2} $ & $-\frac{n+3}{2} - (-2)$ \\
$F_{001} $ & $\infty$ \\
\hline
\end{tabular}

Above, $\mu_0$ is the leading term in the polyhomogeneous conormal expansion of $H_0$ at $F_{100}$ and by symmetry at $F_{010};$  $\mu_0$ is determined by the eigenvalues of the induced Laplacian on $Y,$ the dimension, and the rank of the bundle (see \cite{chee0}).


\subsection{Ac scattering heat kernel}

A summary of the ac scattering heat kernel, space and calculus is given here; for the details, see the appendix.  Let $\bar{Z}$ be a compactified ac scattering space with boundary defined by $\{x=0\}$ and local coordinates $(x,y)$ near the boundary.  Let $\Delta_z$ be the Friedrich's extension of a geometric Laplacian on $Z.$    

\subsubsection{The ac scattering heat space}

First, we construct the ac scattering double space,
$$\bar{Z}^2 _{sc} := \left[ [\bar{Z} \times \bar{Z} ; \partial \bar{Z} \times \partial \bar{Z}] ; \Delta(Y \times Y) \cap F_{110} \right]$$  
where $F_{110}$ is the face created by the first blowup.  This construction comes from \cite{hv}.  Then, the ac scattering heat space is 
$$\bar{Z}^2 _{sc, h} = \left[ \bar{Z}^2 _{sc} \times \R^+ ; \Delta(Z \times Z) \times \{t=0\}, dt \right].$$

The ac scattering heat space has six boundary faces described in the following table.  

\begin{tabular}{|l|l|l|}
\hline
Face & Geometry of face & Defining function in local coordinates \\
\hline
$F_{220}$ & $N^+ ( \Delta(Y \times Y)) \times \R ^+$ & $\rho_{220} = (x^2 + (x')^2 + |y-y'|^2) ^{\frac{1}{2}}$ \\
$F_{220}$ & $N^+((Y \times Y) - \Delta(Y \times Y)) \times \R^+$ & $\rho_{110} = (x^2 + (x')^2)^{\frac{1}{2}}$ \\
$F_{100}$ & $Z \times Y \times \R^+$ & $\rho_{100} = x$ \\
$F_{010}$ & $Y \times Z \times \R^+$ & $\rho_{010} = x'$ \\
$F_{d2}$ & $PN^+ _t( \Delta(Z \times Z))$ & $\rho_{d2} = (|z-z'|^4 + t^2)^{\frac{1}{2}}$ \\
$F_{001}$ & $(Z \times Z) - \Delta(Z \times Z)$ & $\rho_{001} = t$ \\
\hline
\end{tabular}

\subsubsection{Ac scattering heat calculus}

Elements of the ac scattering heat calculus are distributional section half densities of $Z^2 _+$ which are smooth on the interior and lift to be polyhomogeneous on $\bar{Z}^2 _{sc, h}.$  Let $\mu$ be a smooth, non-vanishing half density on $\bar{Z}^2 _+$ and let $\nu$ be a smooth, non-vanishing half density on $\bar{Z}^2 _{sc,h}$.  


\begin{acheatcalc}
For any $k \in \mathbb{R}$ and index sets $E_{110},$ $E_{220},$ $A \in \Psi^{E_{110}, E_{220},  k} _{sc, H} $ if the following hold.

\begin{enumerate}
\item $A \in \mathcal{A}^{-\frac{1}{2} + E_{110}} _{phg}$ at $F_{110} .$ 
\item $A \in \mathcal{A}^{-\frac{n+2}{2} + E_{220}} _{phg}$ at $F_{220} .$ 
\item $A$ vanishes to infinite order at $F_{001},$ $F_{100},$ and $F_{010}.$
\item $A \in \rho_{d2} ^{- \frac{n+3}{2} -k} \mathcal{C}^{\infty} (F_{d2} ).$
\end{enumerate}
\end{acheatcalc}

Two elements of the ac scattering heat calculus compose as follows.   

\begin{accomposition}
Let $A \in \Psi^{A_{110} , A_{220},  k_a} _{sc,H}$, and $B \in \Psi^{B_{110}, B_{220}, k_b} _{sc, H}$.

Then, the composition $B \circ A$ is an element of $\Psi^{A_{110} + B_{110}, A_{220} + B_{220}, k_a + k_b} _{sc, H}$.
\end{accomposition}  
 
The ac scattering heat kernel is constructed analogously to the $b$ and conic heat kernels.  The model heat kernel in this case is the lift of the Euclidean heat kernel to $\bar{Z}^2 _{sc, h}$ and by construction the leading orders of the ac scattering heat kernel are that of the model kernel at the boundary faces of $\bar{Z}^2 _{sc, h}.$  This is stated in the following theorem.  
 
\begin{acparametrix}
Let $(Z, g_z)$ be an asymptotically conic manifold with cross section $(Y, h)$ at infinity.  Let $(E, \nabla)$ be a Hermitian vector bundle over $(Z, g_z)$ which induces a compatible bundle over $(Y,h).$  Let $\Delta$ be a geometric Laplacian on $(Z, g_z)$ associated to the bundle $(E, \nabla).$   Then there exists $H \in \Psi^{E_{110}, E_{220} , -2} _{sc, H} $ satisfying:  

$$ (\partial_t + \Delta) H(z, z', t) = 0,\,  t>0, $$
$$ H(z, z', 0) = \delta (z-z').$$
Moreover, $H$ vanishes to infinite order at $F_{110}$ and is smooth up to $F_{220}.$  
\end{acparametrix}  

The proof of this theorem is in the appendix.  


\section{Heat Kernel Convergence} 

The interaction of the heat kernels will be studied on the asymptotically conic convergence (acc) heat space.  

\subsection{The acc heat space} 
The acc heat space construction is similar to the heat space constructions of section 6 and the double space construction in section 3.  First we let,
$$\cS^2 _h := \left[ [ \cS^2 \times \R^+; C_1 \times C_1 \times \{0\}, dt]; C_1 \times C_1 \right].$$
Let 
$$\cH_0 = \{p \in \cS^2 _h : \epsilon(p) = \epsilon'( p) \},$$
and let $\Delta^* (\cS^2)$ be the lift of the diagonal in $\cS^2 - (C_1 \times C_1)$ to $\cH_0.$  An argument analogous to the smoothness of the acc double space shows that $\cH_0$ is a smooth submanifold of $\cS^2 _h.$  Then, the acc heat space $\cH$ is
$$\cH := [ \cH_0 ; \Delta^*(\cS^2) \times \{0 \}, dt].$$

The acc heat space has several boundary faces of codimension one, however, we are only interested in those at $\epsilon =0.$  These are the following.

\begin{tabular}{|l|l|l|}

\hline
Arising from & face & geometry \\
\hline
$x=0, x'=0,$ & $F_{1010} $ & $\left[ [ [ \bar{Z} \times \bar{Z} \times \R^+ ; Y \times Y \times \{0\}, dt]; Y \times Y ]; \Delta( Z \times Z) \times \{0\}, dt \right]$ \\
$x=0, r'=0$ & $F_{1001} $ & $\left[ [\bar{Z} \times M_0 \times \R^+; Y \times Y \times \{0\}, dt]; Y \times Y  \right]$ \\
$r=0, x'=0$ & $F_{0110}$ & $\left[ [M_0 \times \bar{Z} \times \R^+; Y \times Y \times \{0\}, dt]; Y \times Y  \right]$ \\
$r=0, r'=0$ & $F_{0101}$ & $\left[ [ [ M_0 \times M_0 \times \R^+ ; Y \times Y \times \{0\}, dt]; Y \times Y ]; \Delta(M_0 ^0 \times M_0 ^0) \times \{0 \}, dt \right]$ \\
\hline
\end{tabular}

Note that $F_{1010} \cong \widetilde{Z}^2 _{b,h}$ the $b$-heat space associated to $\bar{Z}$ with an additional benign blowup at the intersection of the boundary in each copy of $\bar{Z}$ at $t=0.$  Similarly, $F_{0101} \cong \widetilde{M}^2 _{0,h}$ the conic heat space with an additional benign blowup at the intersection of the boundary for $t >0.$  These are the faces in which we are most interested and should be thought of as the ``front faces'' while the other two faces should be thought of as the ``side faces.''  The $\epsilon =0$ boundary faces meet in codimension 2 corners described in the following table. 

\begin{tabular}{|l|l|l|l|}
\hline
Arising from & corner & geometry & in faces \\
\hline
$x=0, r=0, x'=0$ & $C_{1110}$ & $[ Y \times \bar{Z} \times \R^+; Y \times Y \times \{0\}, dt]$ & $F_{1010}, F_{0110}.$\\
$x=0, r=0, r'=0$ & $C_{1101}$ & $[Y \times M_0 \times \R^+; Y \times Y \times \{0\}, dt]$ & $F_{1001}, F_{0101}$ \\
$x=0, x'=0, r'=0$ & $C_{1011}$ & $[\bar{Z} \times Y \times \R^+; Y \times Y \times \{0 \}, dt]$ & $F_{1010}, F_{1001}$\\
$r=0, x'=0, r'=0$ & $C_{0111}$ & $[ M_0 \times Y \times \R^+; Y \times Y \times \{0\}, dt]$ & $F_{0101}, F_{0110}.$ \\
\hline
\end{tabular}

There are also higher codimension corners contained in the boundary faces

\subsection{The acc heat calculus}

The acc heat calculus is a parameter ($\epsilon$) dependendent operator calculus which incorporates the smooth, conic, ac scattering and $b$-heat calculi.  Let $\Psi^{k} _{\epsilon, H}$ be the heat calculus of order $k$ for the smooth compact manifold $(M, g_{\epsilon}),$ (see \cite{tapsit} or \cite{moi}). The acc heat calculus consists of kernels that restrict for each $\epsilon$ to an element of $\Psi^{k} _{\epsilon, H}$ and which have a polyhomogeneous expansion at $\epsilon=0$ in terms of elements of $\Psi^{k, E_{110}} _{b,H}$ and $\Psi^{k, E_{112}, E_{100}, E_{010}} _{0, H}.$ 

\begin{accheatcalc}
The \em asymptotically conic convergence heat calculus of order $k,$ \em written $\Psi^{k, E_{1010}, E_{0101}, E_{1001}, E_{0110}} _{acc, H}$ consists of kernels $A$ such that the following hold. 

\begin{enumerate}
\item For each $\epsilon >0,$ $A$ restricts to an element of $\Psi^k _{\epsilon, H}.$
\item In a neighborhood of $F_{1010},$ $A$ has an asymptotic expansion in $\rho_{1010}$ with index set $E_{1010}$ such that the coefficients are elements of the $b$-heat calculus of order $k.$  Such an expansion is of the form
$$A \sim \sum_{j \geq 1} \sum_{0 \leq p_0 \leq p \leq p_j} (\rho_{1010})^{\alpha_j} (\log \rho_{1010} )^{p} A_{j,l} ,$$ 
with $A_{j,l} \in \Psi^{k, E_{110} ^{j}} _{b,H}.$  Above, if for some $j,$ $p_j =0,$ then there are no $\log$ terms. 
\item In a neighborhood of $F_{0101},$ $A$ has an asymptotic expansion in $\rho_{0101}$ with index set $E_{0101}$ such that the coefficients are elements of the conic heat calculus of order $k.$  Such an expansion is of the form
$$A \sim \sum_{l \geq 1}  \sum_{0 \leq p_0 \leq p \leq p_l} (\rho_{0101})^{\alpha_l} (\log \rho_{0101} )^{p} B_{l,p} ,$$ 
with $B_{l,p} \in \Psi^{k, E_{112} ^{l}, E_{100} ^{l}, E_{010} ^{l}} _{0,H}.$
\item $A$ has asymptotic expansion in $\rho_{1001}$ at $F_{1001}$ with index set $E_{1001}$ and asymptotic expansion in $\rho_{0110}$ at $F_{0110}$ with index set $E_{0110}.$
\end{enumerate}
\end{accheatcalc} 

\subsubsection{Composition} 

Let $A \in \Psi^{k_a, A_{1010} ,A_{0101} , A_{1001}, A_{0110}} _{acc, H},$ $B \in \Psi^{k_b, B_{1010} , B_{0101} , B_{1001}, B_{0110}} _{acc, H}.$  We expect $A$ and $B$ to compose as follows.  

For each $\epsilon >0,$ the restrictions of $A$ and $B$ to the $\epsilon$ slice of $\mathcal{H}$ compose to give 
$$(A \circ B) |_{\epsilon} \in \Psi^{k_a + k_b} _{\epsilon, H}.$$
At $F_{1010},$ $A$ has expansion
$$A \sim \sum_{j \geq 1}  \sum_{p_0 \leq p \leq p_j} (\rho_{1010})^{\alpha_j} (\log \rho_{1010} )^{p} A_{j,p} ,$$ 
with $A_{j,p} \in \Psi^{k_a, A_{110} ^{j}} _{b,H}$ and $\{ \alpha_j, p_j \} = A_{1010}.$  Similarly, $B$ has expansion
$$B \sim \sum_{l \geq 1}  \sum_{q_0 \leq q \leq q_l} (\rho_{1010})^{\beta_l} (\log \rho_{1010} )^{q} B_{l,q} ,$$ 
with $B_{l,q} \in \Psi^{k_b, B_{110} ^{l}} _{b,H}$ and $\{ \beta_l , q_l \} = B_{1010},$ so at $F_{1010}$ the composition $A \circ B$ has expansion 
$$A \circ B \sim \sum_{m \geq 1} \sum_{j+l = m}  \sum_{r_0 \leq r \leq r_m} (\rho_{1010})^{\alpha_j + \beta_l} (\log (\rho_{1010}))^{r} A_{j,p} \circ B_{l,q} .$$
By the composition rule for the $b$-heat calculus, $A_{j,p} \circ B_{l,q} \in \Psi^{k_a + k_b, A_{110} ^{j} + B_{110} ^{l}} _{b, H}$ so $A \circ B$ has expansion
$$A \circ B \sim \sum_{j,l \geq 1}  \sum_{r_0 \leq r \leq r_{k}}(\rho_{1010})^{\gamma_{k}} (\log (\rho_{1010}))^{r} C_{k, r}$$
with index set $\{ \gamma_{k} , r_k \} = A_{1010} + B_{1010}$ and $C_{k,r} \in \Psi^{k_a + k_b, A_{110} ^{j} + B_{110} ^{l}} _{b, H}.$

At $F_{0101},$ $A$ has expansion
$$A \sim \sum_{j \geq1}  \sum_{r_0 \leq r \leq r_j}(\rho_{0101})^{a_j} (\log(\rho_{0101}))^{r} R_{j,r}$$
with index set $\{ a_j, r_j \} = A_{0101}$ and $R_{j,r} \in \Psi^{k_a, A_{112} ^j, A_{100} ^j, A_{010} ^j} _{0, H}.$  Similarly at $F_{0101},$ $B$ has expansion
$$B \sim \sum_{l \geq1}  \sum_{s_0 \leq s \leq s_l} (\rho_{0101})^{b_l} (\log(\rho_{0101}))^{s} S_{l,s}$$
with index set $\{ b_l , s_l \} = B_{0101}$ and $S_{l,s} \in \Psi^{k_b, B_{112} ^l, B_{100} ^l, B_{010} ^l} _{0, H}.$  The composition $A \circ B$ then has an expansion at $F_{0101}$ 
$$A \circ B \sim \sum_{m \geq 1} \sum_{j+l = m}  \sum_{p_0 \leq p \leq p_j} (\rho_{0101})^{a_j + b_l} (\log(\rho_{0101}))^{p} R_{j,r} \circ S_{l,s}.$$
By the composition rule for the conic heat calculus $R_{j,r} \circ S_{l,s} \in \Psi^{k_a + k_b, A_{112} ^j + B_{112} ^l, A_{110} ^j + B_{100} ^l, A_{010} ^j + B_{010} ^l} _{0,H}.$  This shows that the composition $A \circ B$ has an expansion of the form
$$A \circ B \sim \sum_{j,l \geq 1}  \sum_{p_0 \leq p \leq p_{j,l} } (\rho_{0101})^{c_{j,l}} (\log(\rho_{0101}))^{p} G_{j,l}$$
with index set $\{ c_{j,l}, p_{j,l} \} = A_{0101} + B_{0101}$ and $G_{j,l} \in \Psi^{k_a + k_b, A_{112} ^j + B_{112} ^l, A_{110} ^j + B_{100} ^l, A_{010} ^j + B_{010} ^l} _{0,H}.$

At $F_{1001}$ and $F_{0110},$ $A$ and $B$ have expansions with index sets $A_{1001}, B_{1001}$ and $A_{0110}, B_{0110},$ respectively.  The composition then has expansions at  $F_{1001}$ and $F_{0110}$ with index sets $A_{1001} + B_{1001}$ and $A_{0110} + B_{0110},$ respectively.

The composition rule is not required for the proof of our main theorem, but we expect it to follow as above.  A sketch of the proof including construction of the acc triple heat space, the key technical tool for proving the composition rule, is included in appendix B.  


\begin{bigthm2}
Let $(M_0, g_0 )$ be a compact Riemannian $n$-manifold with isolated conic singularity,  and let $(Z, g_z)$ be an asymptotically conic space, with $n \geq 3.$  Assume $(M, g_{\epsilon} )$ converges asymptotically conically to $(M_0 , g_0).$  Let $(E_0, \nabla_0)$ and $(E_z, \nabla_z)$ be Hermitian vector bundles over $(M_0, g_0)$ and $(Z, g_z),$ respectively, so that each of these bundles in a neighborhood of the boundary is the pullback from a bundle over the cross section $(Y, h).$  Let $\Delta_0,$  $\Delta_z$ be the corresponding Friedrich's extensions of geometric Laplacians, and let $\Delta_{\epsilon}$ be the induced geometric Laplacian on $(M, g_{\epsilon}).$  Then the associated heat kernels $H_{\epsilon}$ have a full polyhomogeneous expansion as $\epsilon \to 0$ on the asymptotically conic convergence (acc) heat space with the following leading terms at the $\epsilon =0$ faces:
$$H_{\epsilon} (z, z', t) \to H_0 (z, z', t')   \textrm{ at } F_{0101} $$
$$H_{\epsilon} (z, z', t) \to (\rho_{1010})^2 H_b (z, z', \tau)  \textrm{ at } F_{1010},$$
$$H_{\epsilon} (z, z', t) = O\left( (\rho_{1001})^2 (\rho_{0110})^2 \right) \textrm{ at } F_{1001}, \, F_{0110}.$$

Above, $H_0$ is the heat kernel for $\Delta_0,$ $H_b$ is the $b$-heat kernel for $\Delta_z$ and $t'$ and $\tau$ are rescaled time variables with 
$$t' = \frac{t}{(\rho_{1010})^2}, \, \tau = \frac{t}{(\rho_{1001} \rho_{0110})^2}.$$   This convergence is uniform in $\epsilon$ for all time and moreover, the error term is bounded by $C_N \epsilon t^N,$ for any $N \in \N_0.$    
\end{bigthm2}

\subsection{Proof}

This proof is modeled after the parametrix construction of \cite{tapsit}.  The proof is in four steps.  First, we lift the operator $\partial_t + \Delta_{\epsilon}$ to $\mathcal{H}.$  Second, we construct the \em acc model heat kernel \em as an element of the acc heat calculus to be a solution kernel for the lifted operator.  Next, we estimate the error term of this solution kernel.  Finally, we use the $b$-heat calculus, introduce the acc conic triple heat space, and use the conic heat calculus to solve away the error term.  

\subsubsection{Lifted heat operator}
To analyze the behavior of the heat kernels under ac convergence we must lift the heat operator $\partial_t + \Delta_{\epsilon}$ depending on the parameter $\epsilon$ to the acc heat space.  Specifically, we are interested in the leading orders of $\partial_t + \Delta_{\epsilon}$ near the $F_{1010}$ and $F_{0101}$ (front) faces.  In a neighborhood of $F_{1010},$ $x$ and $x'$ are small and approaching zero.  So, we may be initially tempted to use the local cordinates $(r, y, r', y', t, \epsilon).$  However, these are not good coordinates up to the face $F_{1010}$ since $r, r', \epsilon$ vanish at various regions of this face.  So instead we consider the projective coordinates
$$\sigma = \frac{r}{r'}, \, \sigma' = r', \, \eta = \frac{\epsilon}{\sigma \sigma'}, \, \tau = \frac{t}{\eta^2}.$$
Near $F_{1010}$ in the coordinates $(r, y, r', y', t, \epsilon)$ the heat operator has the form
$$\partial_t + \epsilon^{-2} ( r^4 \partial_r ^2 + r^2 \Delta_y + lot)$$
where $lot$ are lower order terms, in the sense that they contain spatial derivatives of lower order.  In projective coordinates this becomes
$$\eta^{-2} [ \partial_{\tau} + \Delta_{b, \sigma}],$$
where 
$$\Delta_{b, \sigma} = (\sigma \partial_{\sigma} )^2 + \Delta_y + lot,$$
is an elliptic order $2$ $b$-operator on $(0,\infty)_{\sigma} \times Y.$  Note that these projective coordinates are valid up to and on $F_{1010}$ away from the $r'=0$ face and corners meeting this face.  In these projective coordinates $\eta$ is a defining function for the side face $F_{1001}.$  By symmetry, the behavior in a neighborhood of the $\sigma'=0$ corner of $F_{1010}$ is the same behavior as at the $\sigma =0$ corner of $F_{1010}.$

Near $F_{0101}$ in the coordinates $(x, y, x', y', t, \epsilon)$ the heat operator has the form
$$\partial_t + \partial_x ^2 + x^{-2} \Delta_y + lot.$$
Here we consider the projective coordinates
$$s = \frac{x}{x'}, \, s'= x', \, t' = \frac{t}{(s')^2}, \, \xi = \frac{\epsilon}{s'}.$$
The heat operator in this coordinates is
$$(s')^{-2} (\partial_{t'} + \Delta_{0,s}),$$
where $\Delta_{0,s} = \partial_s ^2 + s^{-2} \Delta_y + lot$ is a geometric Laplacian on $(0, \infty)_s \times Y.$  Note that these projective coordinates are valid up to and on $F_{0101}$ away from the $x'=0$ face and corners meeting this face.  Since $s'$ is a defining function for the side face $F_{0110}$ as well as $F_{1010}$ in these coordinates we see that the heat operator lifted to the acc heat space has the following leading order behavior at the boundary faces.  

\begin{tabular}{|l|l|l|}
\hline
Face & $\partial_t + \Delta_{\epsilon}$ leading term & time variable \\
\hline
$F_{1010}$ & $(\rho_{1010})^{-2} (\partial_{\tau} + \Delta_{b, \sigma})$ & $\tau = \frac{t}{(\rho_{1001} \rho_{0110})^2}$ \\
$F_{0101}$ & $\partial_{t'} + \Delta_{0,s}$ & $t' = \frac{t}{(\rho_{1010})^2 }$ \\
$F_{0110}$ & $(\rho_{0110})^{-2}$ & \\
$F_{1001}$ & $(\rho_{1001})^{-2}$ & \\
\hline
\end{tabular} 

We can now construct the acc model heat kernel $H_1$ as an element of the acc heat calculus that solves the heat equation at the $\epsilon =0$ boundary faces (and interior) of $\mathcal{H}$ to at least first order. 

\subsubsection{Acc model heat kernel}

At $F_{0101},$ let $H_1 = H_0 (z, z', t'),$ where $H_0$ is the heat kernel for $(M_0 , g_0)$ with rescaled time variable $t'.$  Extend $H_1$ smoothly off $F_{0101},$ and at $F_{1010}$ let $H_1$ be asymptotic to  
$$H_1 \sim (\rho_{1010})^2 H_b (z, z', \tau),$$
where $H_b$ is the $b$-heat kernel with rescaled time variable $\tau,$ and let $H_1$ restrict to $H_b$ at $F_{1010}.$  At $F_{1001}$ and $F_{0110},$ let $H_1$ be asymptotic to 
$$H_1 \sim (\rho_{1001})^{2} \textrm{ at $F_{1001}$, } \, H_1 \sim (\rho_{0110})^{2} \textrm{ at $F_{0110}. $}$$ 
Extend $H_1$ smoothly off the boundary faces so that for $\epsilon > \epsilon_0 >0,$ $H_1$ is the heat kernel $H_{\epsilon}$ for $(M, g_{\epsilon}).$  As $t \to 0$ away from the diagonal $H_1$ vanishes to infinite order.  At $t=0$ along the diagonal, a purely local construction away from the boundaries and singularities using the Euclidean heat kernel solves the heat equation up to $O(t).$  The standard parametrix construction \cite{tapsit} using the jet of the metric at the base point of each fiber improves this to an error term vanishing to infinite order in $t$ as $t \to 0$ uniformly down to $\epsilon =0.$  Since the setting is a bit different here, we include this construction.  

\subsubsection{Acc model heat kernel construction along diagonal at $t=0$}

In a neighborhood of the faces diffeomorphic to $PN^+(\Delta)$ in $\cH$ we carry out a local construction as in \cite{tapsit} chapter 7.  Let $X$ be a manifold; since this construction is the same for  $X=M,$ $X=M_0$ and $X=\bar{Z}$ we use $X$ to simplify notation.  Let $F_X$ denote the $PN^+(\Delta(X))$ face in $\cH,$ where $\Delta(X)$ is the diagonal in $X \times X.$  Let $N(\partial_t + \Delta_X)$ be the restriction to $F_X$ of the lift of $(\partial_t + \Delta_X)$ to $\cH,$ where $\Delta_X$ is our geometric Laplacian on $X.$  With the heat calculus normalization at $F_X,$ an element $A$ of the acc heat calculus of order $-k$ restricts to $F_X$ as follows
$$N(A) = t^{(k + n + 2)/2} A|_{F_X}.$$

As in \cite{tapsit} we observe that $F_X$ is naturally diffeomorphic to a radial compactification of the tangent space of $X,$ with each fiber of $F_X$ over $(x, x, 0)$ diffeomorphic to the tangent space at $x,$ $fiber(x) \cong T_x X.$  Then, from \cite{tapsit} 7.15,  
\begin{equation} \label{normal}
N(t (\partial_t + \Delta_X) A) = [ \sigma(\Delta_X) - \frac{1}{2} (R + n + k + 2)] N(A),\end{equation}
where $R$ is the radial vector field on the fibers of $TX.$  Note that if $G_0$ satisfies 
$$t(\partial_t + \Delta_X)G_0 = O(t^{\infty}) \textrm{ as } t \to 0, \quad G_0 |{t=0} = \delta(x - x'),$$
then $G_0$ also satisfies
\begin{equation}\label{normalcondit}(\partial_t + \Delta_X)G_0 = O(t^{\infty}) \textrm{ as } t \to 0, \quad G_0 |{t=0} = \delta(x - x'). \end{equation}
So, we may work with $t(\partial_t + \Delta_X)$ as in \cite{tapsit}.  Our initial parametrix $G_0$ will have order $k=-2$ at $F_X.$  Then, from \ref{normal} we have the following equation for $G_0$
$$[\sigma(\Delta_X) - \frac{1}{2}(R + n)] N(G_0) = 0.$$
From \cite{tapsit} 7.13, in order for $G_0$ to satisfy the initial condition it must satisfy,  
$$\int_{\textrm{fiber}} N(G_0) =1.$$
Since these conditions are fiber-by-fiber we introduce local coordinates so that 
$$\sigma(\Delta_X) = D_1 ^2 + \ldots + D_n ^2 \, \, \textrm{ on } T_x X.$$
Then we have
\begin{equation} \label{localcondit} \left[ D_1 ^2 + \ldots + D_n ^2 - \frac{1}{2}(R + n) \right] N(G_0) =0,\end{equation}
so 
\begin{equation} \label{g0} N(G_0) = (2\pi)^{-\frac{n}{2}} \exp \left(- \frac{|X|_x ^2}{4}\right) \end{equation}
is the desired solution, where $X$ is a projective local coordinate on $F_x,$ $X = \frac{x-x'}{t^{1/2}}$ (see \cite{tapsit} (7.36)), and $| * |_x$ is the Riemannian norm on $TX$ induced by the metric at $x.$  To see that this is the desired solution, consider the Fourier transform of (\ref{localcondit}) with $u = N(G_0) |_{T_x X},$  
$$(\xi \partial_{\xi} + 2 |\xi|^2) \hat{u} = 0, \quad \hat{u}(0) = 1.$$
Then by standard results in ordinary differential equations, the expression in (\ref{g0}) is the unique decaying solution.  

Now we may iterate this to solve up to higher order.  Assume we have found $G_0, \ldots, G_k$ satisfying 
$$t(\partial_t + \Delta_X) G_j = R_j ,$$
where $R_j$ is of order $-3 - j$ at $F_X.$  To find $G_{k+1} = G_k - T_k,$ we wish to solve
$$t (\partial_t + \Delta_X) T_k = R_k + R_{k+1},$$
where we have already found $R_k$ of order $-3-k,$ and $R_{k+1}$ will be of order $-4-k.$ Lifting to $TX$ this becomes 
$$\left[ \sigma(\Delta_X) - \frac{1}{2}(R + n -j -1) \right] N(T_k) = N(R_k),$$
which we may again solve via Fourier transform.  Letting $u = N(T_k)$ and $f = N(R_k)$ we find 
$$\hat{u} (\xi) = \int_0 ^1 \exp((r-1)|\xi|^2)\hat{f} (r \xi) r^{k+1} dr$$
is the desired solution.  This completes the inductive construction for all $k.$  Now the successive $T_j = G_{j+1} - G_j$ give a formal power series at $F_X$ which can be summed by Borel's Lemma so that $G$ is order $-2$ at $F_X$ and satisfies (\ref{normalcondit}).  We then set the acc model heat kernel $H_1 = G$ in a neighborhood of $F_X.$  Since this construction is the same for $X = M,$ $X = M_0$ and $X=\bar{Z}$ and for each the error term vanishes to infinite order as $t \to 0,$ the error term for $H_1$ has infinite order vanishing as $t \to 0$ for all $\epsilon \geq 0.$  

\subsubsection{Error term approximation}
For each $\epsilon >0,$ let $E(z, z', t, \epsilon) = H_{\epsilon} (z, z', t) - H_1 (z, z', t, \epsilon).$  Let $K$ be defined for each $\epsilon >0$ by 
$$(\partial_t + \Delta_{\epsilon}) E (z, z', t, \epsilon) = K(z, z', t, \epsilon).$$
By construction of $H_1,$ $K = O(\epsilon t^{\infty})$ as $\epsilon, t \to 0,$ so for any$N \in \mathbb{N},$ there is $C >0$ such that for any $(z, z') \in M \times M,$ 
$$|K(z, z', t, \epsilon)| < C \epsilon t^N.$$  
Moreover, $K$ has a polyhomogeneous expansion down to $\epsilon =0.$

For each $\epsilon >0,$ $E$ is smooth on $\mathcal{H}$ for $t  >0$ by parabolic regularity applied for each $\epsilon >0$ since $K$ is $O(t^{\infty}).$  By construction $E$ is smooth down to $t=0,$ so $E(z, z', t, \epsilon)$ is smooth on the blown down space, $M \times M \times \mathbb{R}^+ \times (0,\delta]_{\epsilon}.$  The following maximum principle argument on $M \times [0,T]_t$ shows that $E$ is also $O(\epsilon t^{\infty})$ as $\epsilon, t \to 0$ in the same sense as $K.$  

Fix $\epsilon >0,$ $z' \in M.$  Since $K = O(\epsilon t^{\infty}),$ fix $C >1$ and $\mathbb{N} \ni N >>1$ such that $|K(z, z', t, \epsilon)|^2 \leq C \epsilon^2 t^{2N}$ for all $z \in M.$  Let $u(z, t)$ = $|E(z, z', t, \epsilon)|^2.$  Let $\Delta$ be the scalar Laplacian for $(M, g_{\epsilon}).$  Then $u$ satisfies
$$(\partial_t + \Delta)u = 2\langle (\partial_t + \nabla^* \nabla)E, E \rangle - | \nabla E|^2 \, = \, 2\langle K - \mathcal{R}E, E \rangle - |\nabla E|^2 $$
$$\leq 2 \langle K, E \rangle \leq 2|K| |E| \leq |K|^2 + |E|^2 = |K|^2 + u.$$
Above we have used the positivity of $\mathcal{R}$ and the compatibility of the bundle connection with the metric.  Now, let $\tilde{u} = e^{-t}u.$  Then $\tilde{u}$ satisfies
$$(\partial_t + \Delta) \tilde{u} \leq e^{-t} |K|^2 \leq C \epsilon^2 t^{2N }.$$
Let $w = \tilde{u} - C \epsilon^2 t^{2N+1}.$  Since $E$ and hence $u$ and $\tilde{u}$ vanish at $t=0,$ $w|_{t=0} = 0$ and $w$ satisfies
$$(\partial_t + \Delta)w \leq C \epsilon^2 t^{2N} - C(2N+1) \epsilon^2 t^{2N} < 0.$$
Fix $T>0$ and consider $w$ on $M \times [0,T]_t.$  If $w$ has a local maximum for $z \in M$ and $t \in (0, T)$ then 
$$(\partial_t + \Delta)w >0,$$
and this is a contradiction.  If $w$ has a maximum at $t=T$ then $\partial_t w \geq 0$ and 
$$(\partial_t + \Delta)w >0,$$
which is again a contradiction.  Therefore, the maximum of $w$ occurs at $t=0$ and so 
$$w \leq C \epsilon^2 t^{2N+1}.$$
This implies 
$$u \leq e^T C \epsilon^2 t^{2N}, \quad \textrm{for } 0< t \leq T,$$
which in turn implies that $E = O(\epsilon t^{N})$ as $\epsilon, t \to 0,$ for any $N \in \mathbb{N}.$ 

\subsubsection{Full Construction} 
We solve away the error term of the acc model heat kernel $H_1$ and construct the full asymptotic expansion of the acc heat kernel $H_{acc}$ on $\mathcal{H}.$  First consider $F_{1010} \cong \widetilde{Z}^2 _{b,h}.$  This face has the following geometry.

\begin{tabular}{|l|l|l|}
\hline
Boundary Face & Geometry of Face & Arising from \\
\hline
$F_{bd2}$ & $PN^+ ( \Delta (Z \times Z) )$ & parabolic blowup of diagonal at $t=0$\\ 
$F_{b112}$ & $PN^+(Y \times Y)$ & parabolic blowup of $Y \times Y$ at $t=0$\\
$F_{b110}$ & $SN^+(Y \times Y) \times \R^+$ & blowup of $Y \times Y$ for all time \\
$F_{b100}$ & $Y \times Z \times \R^+$ & boundary in first copy of $Z$\\
$F_{b010}$ & $Z \times Y \times \R^+$ & boundary in second copy of $Z$\\
$F_{b001}$ & $(Z \times Z) - \Delta(Z \times Z)$ & $t=0$ away from diagonal\\
\hline
\end{tabular}

Since at this face $H_1$ is asymptotic to $H_{b} (z, z', \tau),$ at the boundary face where $\tau$ vanishes away from the diagonal, $H_1$ vanishes to infinite order.  Recall 
$$\tau = \frac{t}{(\rho_{1001} \rho_{0110})^2},$$
so $\tau$ vanishes at $F_{b001}.$  Moreover, $H_b$ vanishes to infinite order at the side faces $F_{b100}$ and $F_{b010}.$  At the diagonal face $F_{bd2}$ we've solved $H_1$ up to error vanishing to infinite order in $t.$  So, we have at this point an approximation $H_1$ whose error vanishes to infinite order on the interior of $F_{1010}$ and at all boundary faces except $F_{b112}$ and $F_{b110}.$  Since the face $F_{b112}$ arises from $F_{b110}$ and the blowup of these two faces may be performed in either order,  we solve the error term up to infinite order on the blown-down space in which $F_{b112}$ has been blown down.  Then, lifting the solution to the fully blown up space, the error term will still vanish to infinite order.  The indicial operator for $\Delta_{b, \sigma}$ at $F_{b110}$ is 
$$(\sigma \partial_{\sigma})^2 + \Delta_y$$
on $\R^+ _{\sigma} \times Y.$  We would like to solve 
$$(\partial_{\tau} + \Delta_b) |_{F_{b110}} u = -K|_{F_{b110}},$$
so that $u$ is polyhomogeneous on $\widetilde{Z}^2 _{b,h};$ then $H_1 + u$ would vanish identically at $F_{b110}.$  Since $H_1$ is polyhomogeneous on $F_{1010}$ and smooth up to these boundary faces, the error term $K$ is also.  Expanding $K,$ where we use simply $\rho$ for the projective defining function for $F_{b110},$ 
$$K \sim \sum_{j \geq 0} (\rho )^j k_j, \quad k_j \in \mathcal{C}^{\infty} _0,$$
and expanding the desired solution $u,$
$$u \sim \sum _{j \geq 0} (\rho)^j u_j$$
we may use either separation of variables expanding in eigenfunctions of $\Delta_y$ or the Mellin transform (see \cite{tapsit}) to find $u_0$ satisfying  
$$(\partial_{\tau} + I (\Delta_b)) u_0 = -k_0,$$
with $u_0$ vanishing to infinite order as $\sigma \to 0, \infty:$  at the side faces $F_{b100}, F_{b010}.$    Since $\Delta_b - I(\Delta_b) = (\rho) (L_1),$ where $L_1$ is also a $b$-differential operator, we may now iteratively solve for $u_1, u_2, \ldots,$ to solve the equation to increasingly higher order.  Recall that in the projective coordinates near this face, $\sigma'$ defines $F_{b110}$ and since the operator does not differentiate with respect to $\sigma' = r'$ the defining function commutes past the operator.  Using Borel summation we construct $u$ so that 
$$(\partial_{\tau} + \Delta_b) u = -K + K_2$$
where $K_2$ vanishes to infinite order at $F_{b110}, F_{b112}.$  Using a smooth cutoff function $\chi$ 
supported in a neighborhood of these faces, the second approximation $H_2 = H_1 + \chi u$ now satisfies
$$(\partial_{\tau} + \Delta_b) H_2 |_{F_{11}} = K_2$$
where $K_2$ vanishes to infinite order on both the interior and all boundary faces of $F_{1010}.$

Since the error now vanishes to infinite order at all boundary faces except those arising from $F_{0101} \cong \widetilde{M}^2 _{0,h},$ we restrict attention to this face.  It is convenient to use the conic heat calculus composition rule, but this requires the conic triple space.  Since our error vanishes to infinite order at all boundary faces except $F_{0101},$ we construct a partial acc triple heat space, the acc conic triple heat space, which contains the conic triple heat space so that we may use the conic heat calculus composition rule.  

\subsubsection{Acc conic triple heat space}
The acc conic triple heat space $\mathcal{H}^3 _c$ is a submanifold constructed from $\mathcal{S}^3 \times \R^+ \times \R^+$ by eight blowups.  Let $X, X', X''$ denote the three copies of the submanifold $X$ in $X ^3.$  Then the blowups are listed in the following table in the order which the blowups are performed together with the name of the face created.  

\begin{tabular}{|l|l|}
\hline
Blowup & Face\\ 
\hline
$Y \times Y' \times Y'' \times \{t=0, t'=0\}, dt, dt'$ & $F_{11122}$ \\
$Y \times Y' \times \{t=0\}, dt$ & $F_{11020}$\\
$Y' \times Y'' \times \{t' =0\}, dt'$ & $F_{01102}$\\
$Y \times Y'' \times \{t'' = |t-t'| = 0 \}, dt''$ & $F_{10122}$\\
$\Delta(\mathcal{S} \times \mathcal{S}' \times \mathcal{S}'') \times \{ t=0, t'=0 \}, dt, dt'$ & $F_{d3}$\\
$\Delta(\mathcal{S} \times \mathcal{S} ') \times \{ t=0 \}, dt$ & $F_{d20}$\\
$\Delta(\mathcal{S} ' \times \mathcal{S} '') \times \{t' =0 \}, dt'$ & $F_{d02}$\\
$\Delta(\mathcal{S} \times \mathcal{S} '' ) \times \{t'' =0 \}, dt''$ & $F_{d22}$\\
\hline
\end{tabular}

Let $\beta^* H_2$ be the lift of $H_2$ to $\cH^3 _c,$ and let $\beta^* K_2$ be the lift of $K_2$ to $\cH ^3 _c.$  Then, $\beta^* K_2$ vanishes to infinite order at all boundary faces except those arising from the lift of $F_{0101}.$  Now let 
$$H_3 := \beta_* (\beta^* H_2 - \beta^* H_2 \beta^* K_2),$$
where $\beta_*$ is the push forward to $\cH$ from $\cH ^3 _c.$  Since $\beta^* K_2$ vanishes to infinite order at all boundary faces except $F_{0101},$ the push forward of $(\beta^* H_2 )(\beta^* K_2)$ to $\cH$ vanishes to infinite order at all boundary faces except $F_{0101},$ where the result is given by the conic heat calculus composition rule.  Consequently, 
$$H_3 = H_2 - \beta_* (\beta^* H_2 \beta^* K_2)$$
vanishes to higher order at the boundary faces of $F_{0101},$ by the conic heat calculus composition rule.  Continuing this construction and using Borel summation, we arrive at $H_{\infty} $ with expansion asymptotic to $H_2, H_3, \ldots$ and satisfying
$$(\partial_{t'} + \Delta_0) H_{\infty} = K_{\infty},$$
where $K_{\infty}$ now vanishes to infinite order on $F_{0101}.$  Using a smooth cutoff function we now have $H_{\infty}$ defined on all of $\mathcal{H}$ satisfying
$$(\partial_t + \Delta_{\epsilon}) H_{\infty} = K_{\infty}$$
where $K_{\infty}$ vanishes to infinite order at all boundary faces of $\mathcal{H}.$  

\subsubsection{Solving away the residual error term}

To complete this construction we must remove the residual error term which vanishes to infinite order at the boundary faces of $\cH.$  It is now convenient to consider the elements of the acc heat calculus as $t$-convolution operators acting on $\cS \times \R^+.$  For an element $A$ which vanishes to infinite order at the boundary faces of $\cH$ and $u$ a smooth half density section of $\cS \times \R^+,$ the $t$-action of $A$ on $u$ is
\begin{equation} \label{Au} Au(t) = \int_0 ^t \langle Au(t-s), u(s) \rangle ds, \end{equation}
where the spatial variables have been suppressed.  As a function of $s', s \geq 0,$ $[Au(s')](s)$ vanishes to infinite order at $s'=0.$  Restricting to $s'=t-s,$
$$[Au(t-s)](s) = s^{-k/2 -1} (t-s)^{j} u_{k,j} (t-s, s),$$
for any $-k,j \in \N_0$ with $u_{k,j}$ a smooth half density, so for any $-k \geq 1$ this is integrable and consequently, $Au(t)$ as in (\ref{Au}) is smooth in $t$ and vanishes rapidly as $t \to 0.$  So, an element $A$ of the acc heat calculus which vanishes to infinite order at all boundary faces of $\cH$ gives rise to a Volterra operator.  Since as a $t$-convolution operator we have
$$(\partial_t + \Delta) H_{\infty} = Id - K_{\infty},$$
we would like to invert $(Id - K_{\infty}).$  Formally, the inverse should be
$$(Id - K_{\infty})^{-1} = \sum_{j \geq 0} K_{\infty} ^j,$$
where $K_{\infty} ^j$ is the $j$-fold composition of $K_{\infty}.$  To show that this Neumann series converges, we estimate the kernel of $K_{\infty} ^j.$  Since $K_{\infty}$ vanishes to infinite order at all boundary faces of $\cH$ we may restrict to submanifolds of $\cH,$ estimating as in \cite{tapsit} and then combine these estimates to estimate $K_{\infty} ^j$ on $\cH.$  The kernel $k_{\infty} ^j$ of the restriction of $K_{\infty} ^j$ to $M \times M \times \R^+ \times \{ \epsilon \}$ is bounded by 
$$|k_{\infty} ^j (z, z', t, \epsilon)| \leq C_{\epsilon, j} \frac{t^{j}}{(j+1)!}, \quad t < T.$$
This follows from the composition rule for the heat calculus on $M$ and the analogous bound in \cite{tapsit} 7.3, where we have taken $k=-2,$ with $k$ as above which we are free to choose since $A$ vanishes to infinite order.  Similarly, by the composition rule for the heat calculus on $M_0$ and the same estimate of \cite{tapsit}, the kernel of the restriction of $K_{\infty} ^j$ to $M_0 \times M_0 \times \R^+$ is bounded by
$$|k_{\infty} ^j (z, z', t) |_{F_{0101}} \leq C_{0, j} \frac{t^{j}}{(j+1)!}, \quad t < T.$$
Similarly, the kernel of the restriction of $K_{\infty} ^j$ to $\bar{Z} \times \bar{Z} \times \R^+$ is bounded by 
$$|k_{\infty} ^j (z, z', t) |_{F_{1010}} \leq C_{z, j} \frac{t^{j}}{(j+1)!}, \quad t < T.$$
These three bounds imply that the constants $C_{\epsilon , j}$ stay bounded as $\epsilon \to 0$ and so we have the following global bound for the kernel of $K_{\infty} ^j$ on both $\cH$ and the blown-down space $\{ \epsilon = \epsilon' \} \subset \cS \times \cS \times \R^+$
$$| k_{\infty} ^j (z, z', t, \epsilon )| \leq C_{j} \frac{t^{j}}{(j+1)!}, \quad t < T.$$
It follows that the Neumann series for $(Id - K_{\infty})^{-1}$ is summable and has an inverse which as a $t$-convolution operator is also of the form $(Id - A)$ where $A$ is an element of the acc heat calculus that vanishes to infinite order at the boundary faces of $\cH.$  Then, the full acc heat kernel is
$$H = H_{\infty} (Id - K_{\infty})^{-1}.$$
As a consequence of this construction, $H$ has a fully polyhomogeneous expansion down to $\epsilon =0$ with leading order terms given by the acc model heat kernel.  
\heart

\textbf{Remarks}

A consequence of this theorem is the convergence 
$$H_{\epsilon} \to H_{0, t'} + \epsilon^2 H_{b, \tau} \quad \textrm{as } \epsilon \to 0,$$
with error term bounded by 
$$C_N \epsilon t^N, \quad \textrm{for any } N \in \mathbb{N}_0, \quad \textrm{and for all $t.$}$$

\begin{appendix}
\section{Asymptotically conic scattering heat kernel}
 
Let $\bar{Z}$ be a compactified ac scattering space with boundary defined by $\{x=0\}$ and local coordinates $(x,y)$ near the boundary.  Let $\Delta_z$ be the Friedrich's extension of a geometric Laplacian on $Z.$  We motivate the definition of the acc heat space by lifting the Euclidean heat kernel to $\bar{Z}^2 _+.$  

Recall the Euclidean heat kernel for $\mathbb{R}^n$,
$$G(z, z', t) = (4 \pi t)^{-\frac{n}{2}} e^{-\frac{|z-z'|^2}{2t}} .$$
Here the coordinate $z = (r,y)$ has not been compactified.  With the compactification of $Z$ given by $x = \frac{1}{r}$ in the local coordinates $(x, y, x', y', t)$ on $\bar{Z}^2 _+$ near the boundary of $\bar{Z}$ the Euclidean heat kernel is
$$G(x, y, x', y', t) = (4 \pi t)^{-\frac{n}{2}} \textrm{ exp}\left(- \frac{ \left( | \frac{1}{x} - \frac{1}{x'} |^2 + |y-y'|^2 \right)}{2t} \right). $$ 
This motivates blowing up
$$S_{110} = \{ (x, y, x', y', t) \, : \, x=0, x'=0 \}.$$ 
In the projective coordinates $s = \frac{x}{x'},$ $s' = x'$ the Euclidean heat kernel is
$$G(s, y, s', y', t) = (4 \pi t)^{-\frac{n}{2}} \textrm{ exp}\left(- \frac{ \left( | \frac{s-1}{ss'} |^2 + |y-y'|^2 \right)}{2t} \right). $$ 
This motivates a second blowup at $s=1$, along the submanifold where the diagonal in $\bar{Z} \times \bar{Z}$ meets the first blown up face
$$S_{220} = \{ (x, y, x', y', t) \, : \, x=0, x'=0, y=y' \}.$$

\subsection{The Ac scattering heat space}

As motivated above, the ac scattering heat space is constructed from $\bar{Z}^2 _+ $ by performing three blowups.  

First, the scattering double space, $\bar{Z}^2 _{sc}$ is constructed:
$$\bar{Z}^2 _{sc} := \left[ [\bar{Z} \times \bar{Z} ; \partial \bar{Z} \times \partial \bar{Z}] ; \Delta(Y \times Y) \cap F_{110} \right]$$
where $F_{110}$ is the face created by the first blowup.  Including the time variable we perform one more blowup to construct the ac scattering heat space,
$$\bar{Z}^2 _{sc, h} = \left[ \bar{Z}^2 _{sc} \times \R^+ ; \Delta(Z \times Z) \times \{t=0\}, dt \right].$$

The ac scattering heat space has six boundary faces described in the following table.  

\begin{tabular}{|l|l|l|}
\hline
Face & Geometry of face & Defining function in local coordinates \\
\hline
$F_{220}$ & $N^+ ( \Delta(Y \times Y)) \times \R ^+$ & $\rho_{220} = (x^2 + (x')^2 + |y-y'|^2) ^{\frac{1}{2}}$ \\
$F_{220}$ & $N^+((Y \times Y) - \Delta(Y \times Y)) \times \R^+$ & $\rho_{110} = (x^2 + (x')^2)^{\frac{1}{2}}$ \\
$F_{100}$ & $Z \times Y \times \R^+$ & $\rho_{100} = x$ \\
$F_{010}$ & $Y \times Z \times \R^+$ & $\rho_{010} = x'$ \\
$F_{d2}$ & $PN^+  _t ( \Delta(Z \times Z))$ & $\rho_{d2} = (|z-z'|^4 + t^2)^{\frac{1}{2}}$ \\
$F_{001}$ & $(Z \times Z) - \Delta(Z \times Z)$ & $\rho_{001} = t$ \\
\hline
\end{tabular}

\subsection{Ac scattering heat calculus}

Elements of the ac scattering heat calculus are distributional section half densities on $Z^2 _+$ which are smooth on the interior and lift to be polyhomogeneous on $\bar{Z}^2 _{sc, h}.$  Let $\mu$ be a smooth, non-vanishing half density on $\bar{Z} \times \bar{Z} \times \mathbb{R}^+$ and let $\nu$ be a smooth, non-vanishing half density on $\bar{Z}^2 _{sc,h}.$  


\begin{acheatcalc}
For any $k \in \mathbb{R}$ and index sets $E_{110},$ $E_{220},$ $A \ in \Psi^{E_{110}, E_{220},  k} _{sc, H} $ if the following hold.

\begin{enumerate}
\item $A \in \mathcal{A}^{-\frac{1}{2} + E_{110}} _{phg}$ at $F_{110} .$ 
\item $A \in \mathcal{A}^{-\frac{n+2}{2} + E_{220}} _{phg}$ at $F_{220} .$ 
\item $A$ vanishes to infinite order at $F_{001},$ $F_{100},$ and $F_{010}.$
\item $A \in \rho_{d2} ^{- \frac{n+3}{2} -k} \mathcal{C}^{\infty} (F_{d2} ).$
\end{enumerate}

\end{acheatcalc}


Elements of the ac scattering heat calculus are Schwartz kernels of operators acting on sections of $Z$ in the usual way and on sections of $Z \times \R^+ _t$ by $t$-convolution.  The composition rule is proven using the ac scattering triple heat space, $\bar{Z}^3 _{sc, h}.$  This space has partial blow down/projection maps to three identical copies of the ac scattering heat space as well as full blow down/projection maps to three identical copies of $\bar{Z}^2 _+;$ these are called the left, right, and center.  Formally, two elements of the ac heat calculus are composed by lifting from the left and right copies of $\bar{Z}^2 _{sc, h}$  to $\bar{Z}^3 _{sc, h},$ multiplying and blowing down/projecting to the center copy of $\bar{Z}^2 _{sc, h}.$  It is key that the triple space be constructed so that these lifts and push-forward maps are $b$-fibrations in order that polyhomogeneity be preserved.  

\subsubsection{The ac scattering triple heat space}
We first construct the ac scattering triple space $\bar{Z}^3 _{sc}$ and later include the time variables.  In a neighborhood of the boundary in each copy of $\bar{Z}$ we have the local coordinates $(x, y),$ which provide the local coordinates $(x, y, x', y', x'', y'')$ on $\bar{Z}^3.$  First we blow up the codimension three corner defined by $\{ x = 0, x' = 0, x'' = 0 \}$.  We call this face $F_{11100}$ with defining function locally given by 
$$\rho_{11100} = (x^2 + (x')^2 + (x'')^2)^{\frac{1}{2}}.$$
Next, we blow up the three codimension two corners corresponding to the $F_{110}$ faces in each of the three copies of $\bar{Z}^2 _{sc,h}.$  These faces are as follows.

\begin{tabular}{|l|l|l|}
\hline
Face & Submanifold to be blown up & Defining Function \\
\hline
$F_{11000}$ & $S_{11000} = \{ x=0, x'=0 \} - F_{11100}$ & $\rho_{11000} = (x^2 + (x')^2)^{\frac{1}{2}}$\\
\hline
$F_{01100}$ & $S_{01100} = \{ x'=0, x''=0 \} - F_{11100}$ & $\rho_{01100} = ( (x')^2 + (x'')^2)^{\frac{1}{2}}$ \\
\hline
$F_{10100}$ & $S_{10100} = \{ x=0, x''=0 \} - F_{11100}$ & $\rho_{10100} = ( (x)^2 + (x'')^2)^{\frac{1}{2}}$ \\
\hline
\end{tabular}

Next we blow up the codimension $2n+1$ corner where the diagonals meet $F_{11100}.$  After the $F_{11100}$ blowup, we have coordinates $( \theta, \theta', \theta'', y, y', y'', \rho_{11100})$, with
$$ x = (\rho_{11100}) \theta, \quad x' =(\rho_{11100}) \theta', \quad x'' =(\rho_{11100}) \theta'', \quad (\theta)^2 + (\theta')^2 + (\theta'')^2 = 1.$$ 
Using these coordinates, we next blow up
$$S_{22200} = \{ \theta = \theta' = \theta'', \, y = y' = y'', \, r_0 = 0 \}.$$
The face created by this blowup is called $F_{22200}$ with defining function
$$\rho_{22200}  \, = \, ( (\theta - \theta')^2 + (\theta' - \theta'')^2 + |y-y'|^2 + |y' - y''|^2 + r_0 ^2 )^{\frac{1}{2}}.$$

After this we blow up the three codimension $n$ corners corresponding to the $F_{220}$ faces in the three copies of the double heat space.  These are as follows.

\begin{tabular}{|l|l|l|}
\hline
Face & Submanifold to be blown up & Defining Function \\
\hline
$F_{22000}$ & $S_{22000} = \{ \theta =0, \theta '=0, y=y' \} $ & $\rho_{22000} = (\theta^2 + (\theta')^2 + |y-y'|^2)^{\frac{1}{2}}$\\
\hline
$F_{02200}$ & $S_{02200} = \{ \theta'=0, \theta''=0, y'=y''\} $ & $\rho_{02200} = ( (\theta')^2 + (\theta'')^2 + |y'-y''|^2)^{\frac{1}{2}}$ \\
\hline
$F_{20200}$ & $S_{20200} = \{ \theta=0, \theta''=0, y=y'' \} $ & $\rho_{20200} = ( (\theta)^2 + (\theta'')^2 + |y-y''|^2)^{\frac{1}{2}}$ \\
\hline
\end{tabular}

We have now constructed the ac scattering triple space, $\bar{Z}^3 _{sc}$.  We next introduce the time variables and perform the parabolic temporal diagonal blowups.  We must first blow up the codimension $2$ corner of $\mathbb{R}^+ \times \mathbb{R}^+$ to preserve symmetry.  Let
$$\mathcal{T}^2 _0 = [ \mathbb{R}^+ \times \mathbb{R}^+ ; \, t = t' = 0 ].$$
The defining function for the blowup of $\{ t = t' = 0 \}$ is $\rho_{00011},$ which we call $t''$ because it plays the role of the third time variable.  We now take $Z^3 _{sc} \times \mathcal{T}^2 _0$ and blow up the temporal diagonal faces.  First, we blow up the codimension $2n+3$ triple diagonal, $S_{d3}$, defined by
$$\{ z = z' = z'', \, t'' = 0 \}.$$
The defining function of this face is $\rho_{d3}$, 
$$\rho_{d3} = ( |z-z'|^4 + |z-z''|^4 + (t'')^2 )^{\frac{1}{4}}.$$
Next, we blow up the three temporal diagonals corresponding to the diagonal faces in the three copies of the double heat space.  These are as follows.

\begin{tabular}{|l|l|l|}
\hline
Face & Submanifold to be blown up & Defining Function \\
\hline
$F_{d20}$ & $S_{d20} = \{ z=z' \} $ & $\rho_{d20} = ( |z-z'|^4 + t^2)^{\frac{1}{4}}$\\
\hline
$F_{d02}$ & $S_{d02} = \{ z'=z'' \}$ & $\rho_{d02} = ( |z'-z''|^4 + (t')^2)^{\frac{1}{4}}$ \\
\hline
$F_{d22}$ & $S_{d22} = \{ z=z'' \}$ & $\rho_{d22} = ( |z-z''|^4 + (t'')^2)^{\frac{1}{4}}$ \\
\hline
\end{tabular}

We have now constructed the ac scattering triple heat space and proceed with the composition rule.  

\begin{accomposition}
Let $A \in \Psi^{A_{110} , A_{220},  k_a} _{sc,H}$, and $B \in \Psi^{B_{110}, B_{220}, k_b} _{sc, H}$.

Then, the composition $B \circ A$ is an element of $\Psi^{A_{110} + B_{110}, A_{220} + B_{220}, k_a + k_b} _{sc, H}$.
\end{accomposition}  

\subsubsection{Proof}
Formally we have,
\begin{equation} \label{eq:alec}
\kappa_{B \circ A} \nu = (\beta_C)_* \left( (\beta_R)^{*} (\kappa_A \nu) (\beta_L)^* (\kappa_B \nu) \right).
\end{equation}

Multiplying both sides of (\ref{eq:alec}) by $\nu$ and using the fact that $(\beta_c )_* (\beta_c )^* (\nu) = \nu$

\begin{equation} \label{eq:aleco}
\kappa_{B \circ A} \nu^2 = (\beta_C)_* \left( (\beta_R)^{*} (\kappa_A \nu) (\beta_L)^* (\kappa_B \nu) (\beta_c )^* (\nu) \right). 
\end{equation}
 
Next we calculate the lifts of the defining functions and half densities from $\bar{Z}^2 _{sc, h}$ to $\bar{Z}^3 _{sc, h}.$  A calculation gives the half density on the heat space $\nu$ in terms of the half density $\mu$ on $\bar{Z}^2 _+$

$$ \nu = (\beta_h)^* \left( (\rho_{110})^{-\frac{1}{2}} (\rho_{220})^{-\frac{n}{2}} (\rho_{d2})^{-\frac{n+1}{2}} \mu \right).$$ 

The ac scattering triple heat space has partial blow down/projection maps $\beta_L,$ $\beta_R,$ and $\beta_C$ to three identical copies of $\bar{Z}^2 _{sc,h}.$  If we denote the three copies of $\bar{Z}$ by $\bar{Z},$ $\bar{Z}',$ $\bar{Z}'',$ and the three time variables $(t, t', t'')$ where $t''$ is from the blowup of $\mathbb{R}^+ \times \mathbb{R}^+$ then the three copies of $\bar{Z}^2 _{sc,h}$ are as follows.  

\begin{tabular}{|l|l|}
\hline
Copy of $\bar{Z}^2 _{sc,h}$ & Associated to in $\bar{Z}^3 _{sc,h}$ \\
\hline
Left & $\bar{Z} \times \bar{Z}' \times \mathbb{R}^+ _t$ \\
\hline
Right & $\bar{Z}' \times \bar{Z}'' \times \mathbb{R}^+ _{t'}$\\
\hline
Center & $\bar{Z} \times \bar{Z}'' \times \mathbb{R}^+ _{t''}$ \\
\hline
\end{tabular}

Next, we compute the lifts of the defining functions for the boundary faces of the heat space to the triple heat space.  

\begin{tabular}{|l|l|l|}
\hline
Lifting map & Defining function on $\bar{Z}^2 _{sc,h}$ & Lift to $\bar{Z}^3 _{sc,h}$ \\
\hline 
$(\beta_L )^*$ & $\rho_{100}$ & $\rho_{10000} \rho_{10100}$ \\
$(\beta_L )^* $& $\rho_{010}$ &$ \rho_{01000} \rho_{01100}$ \\
$(\beta_L )^*$ &$ \rho_{110}$ &$ \rho_{11100} \rho_{11000}$ \\
$(\beta_L )^*$ & $\rho_{220}$ & $\rho_{22200} \rho_{22000}$\\
$(\beta_L )^* $& $\rho_{d2}$ & $\rho_{d3} \rho_{d20}$ \\
$(\beta_L )^*$ & $\rho_{001}$ & $\rho_{00010} \rho_{00011} \rho_{d22} $ \\
\hline
$(\beta_R)^*$ &$ \rho_{100}$ &$ \rho_{01000} \rho_{01100} $\\
$(\beta_R )^*$ &$ \rho_{010}$ &$ \rho_{00100} \rho_{10100}$ \\
$(\beta_R )^*$ & $\rho_{110}$ & $\rho_{11100} \rho_{01100}$ \\
$(\beta_R)^*$ & $\rho_{220}$ & $\rho_{22200} \rho_{02200}$\\
$(\beta_R )^*$ & $\rho_{d2}$ & $\rho_{d3} \rho_{d02}$ \\
$(\beta_R )^*$ & $\rho_{001} $& $\rho_{00001} \rho_{00011} \rho_{d22} $\\
\hline
$(\beta_C )^*$ &$ \rho_{100}$ & $\rho_{10000} \rho_{11000}$ \\
$(\beta_C )^*$ & $\rho_{010}$ & $\rho_{001000} \rho_{01100}$ \\
$(\beta_C )^* $&$ \rho_{110}$ &$\rho_{11100} \rho_{10100}$ \\
$(\beta_C )^*$ & $\rho_{220}$ & $\rho_{22200} \rho_{20200}$\\
$(\beta_C )^* $&$ \rho_{d2}$ & $\rho_{d3} \rho_{d22}$ \\
$(\beta_C )^*$ & $\rho_{001} $& $\rho_{00022} \rho_{00011} \rho_{d22} $ \\
\hline
\end{tabular}

Then,
$$(\beta _L )^* (\nu) = (\beta_L )^* ( (\rho_{110})^{-\frac{1}{2}} (\rho_{220})^{-\frac{n}{2}} (\rho_{d2})^{-\frac{n+1}{2}} \mu ). $$

Next, we use the fact that
$$(\beta_L )^* (\mu) (\beta_R)^* (\mu) (\beta_C )^* (\mu) = \mu_3 ^2 .$$
Here, $\mu^2 _3$ is a smooth density on $\bar{Z} \times \bar{Z} \times \bar{Z} \times \mathbb{R}^+ \times \mathbb{R}^+$, so we may assume
$$ \mu_3 ^2 = \textrm{d}z \textrm{d}z' \textrm{d}z'' \textrm{d}t \textrm{d}t' .  $$

A Jacobian calculation gives the lift of $\mu_3 ^2$ to the triple heat space.  First note
$$(\beta_3 )^* (x) = (\rho_{11100}) (\rho_{11000}) (\rho_{10100}) (\rho_{10000}) ,$$
$$(\beta_3 )^* (x') = (\rho_{11100}) (\rho_{11000}) (\rho_{01100}) (\rho_{01000}) ,$$
$$(\beta_3 )^* (x'') = (\rho_{11100}) (\rho_{01100}) (\rho_{10100}) (\rho_{00100}) .$$

This implies
$$(\beta_3 )^* (\mu_3 ^2 ) = (\rho_{11100})^2 (\rho_{11000} \rho_{01100} \rho_{10100}) ( \rho_{22000} \rho_{02200} \rho_{20200} )^n $$
$$(\rho_{22200} )^{2n+1} (\rho_{d20} \rho_{d02} \rho_{d22})^{n+1} \rho_{d3} ^{2n+3} (t'')  \nu_3 ^2.$$

Here, $\nu_3 ^2$ is a smooth, nonvanishing density on the triple heat space.  Combining this with the above lifts, we arrive at the following formula

$$(\beta_L )^* (\nu) (\beta_R)^* (\nu) (\beta_C )^* (\nu) = (\rho_{11100}) ^{\frac{1}{2}} (\rho_{10100} \rho_{01100} \rho_{10100})^{\frac{1}{2}}$$
$$ (\rho_{22000} \rho_{02200} \rho_{20200})^{\frac{n}{2}} (\rho_{22200}) ^{\frac{n+1}{2}} (\rho_{d3}) ^{\frac{n+3}{2}} (\rho_{d20} \rho_{d02} \rho_{d22} )^{\frac{n+1}{2}} (t'')  \nu^2 _3.$$

To use the push forward theorem of \cite{edge}, we need to write each of these in terms of $b$-densities.  First, we have on the center copy of $\bar{Z}^2 _{sc,h}$
$$ ^b \nu^2 = ( \rho_{100} \rho_{010} \rho_{110} \rho_{220} \rho_{001} \rho_{d2})^{-1} \nu^2 .$$
Then, we have
$$ ^b \nu^2 = (\beta_c)_* (\beta_c)^* ( (\rho_{100} \rho_{010} \rho_{110} \rho_{220} \rho_{001} \rho_{d2} )^{-1}  \nu^2 ).$$
We observe
$$ (\beta_c)^* \left( (\rho_{100} \rho_{010} \rho_{110} \rho_{220} \rho_{001} \rho_{d2})^{-1} \right) = $$
$$(\rho_{10000} \rho_{00100} \rho_{11000} \rho_{01100} \rho_{10100} \rho_{11100} \rho_{22200} \rho_{02200} \rho_{20200} \rho_{d3} \rho_{d22} \rho_{00011}  )^{-1} .$$

So now we multiply both sides of (\ref{eq:aleco}) by 
$(\beta_c)_* (\beta_c)^* (\rho_{100} \rho_{010} \rho_{110} \rho_{220} \rho_{001} \rho_{d2} )^{-1} )$ and inside the right side of (\ref{eq:aleco}) we have 

$$(\rho_{11100} \rho_{11000} \rho{01100} \rho_{10100} )^{-\frac{1}{2}} (\rho_{22000} \rho_{02200} \rho_{22200} )^{\frac{n}{2}} (\rho_{20200}) ^{\frac{n-2}{2}} $$

$$(\rho_{d3}) ^{\frac{n+1}{2}} (\rho_{d20} \rho_{d02} )^{\frac{n+1}{2}} (\rho_{d22}) ^{\frac{n}{2}} (\rho_{10000} \rho_{00100} )^{-1} \nu^2 _3 .  $$
  
To use the push forward theorem, we must change the density $\nu^2 _3$ to a $b$-density.  We observe 
$$ ^b \nu^2 _3 = (  \rho_{11100} \rho_{11000} \rho_{01100} \rho_{10100} \rho_{22200} \rho_{22000} \rho_{02200} \rho_{20200} $$
$$ \rho_{10000} \rho_{01000} \rho_{00100} \rho_{d3} \rho_{d20} \rho_{d02} \rho_{d22} \rho_{00011} \rho_{00010} \rho_{00001} )^{-1} \nu^2_3 .$$

So, we now have for the composition formula 

$$
(\beta_c)_* ( \tilde{\kappa_A } \tilde{\kappa_B }  (\rho_{11100} \rho_{11000} \rho_{01100} \rho_{10100} )^{\frac{1}{2}} (\rho_{22200} \rho_{22000} \rho_{02200} )^{\frac{n+2}{2}}$$
$$(\rho_{20200} ) ^{\frac{n}{2}} (\rho_{d3} \rho_{d20} \rho_{d02} )^{\frac{n+3}{2}} (\rho_{d22}) ^{\frac{n+1}{2}} \rho_{01000} \rho_{00011} \rho_{00010} \rho_{00001} (^b \nu_3 ^2 ) ).
$$

We observe the following orders of $\tilde{\kappa_A}$ on $\bar{Z}^3 _{sc,h}.$

\begin{tabular}{|l|l|}
\hline
Face & $\tilde{\kappa_A}$ Index Set/Leading Order \\
\hline
$F_{11100}$ &  $-\frac{1}{2} + A_{110}$ \\
$F_{11000}$ &  $ -\frac{1}{2} + A_{220}$  \\
$F_{01100},$ $F_{10100},$ $F_{02200},$ $F_{20200},$ $F_{d22}$ & $\infty$ \\
$F_{22200},$ $F_{22000}$ & $-\frac{n+2}{2} + A_{220}$\\
$F_{d3},$ $F_{d20}$ & $-\frac{n+3}{2} - k_a$ \\
$F_{10000},$ $F_{01000},$ $F_{00100},$ $F_{00010},$ $F_{00011}$ &$ \infty$  \\
\hline
\end{tabular}    

Similarly, for $\tilde{\kappa_B}$ we have orders as follows.

\begin{tabular}{|l|l|}
\hline
Face & $\tilde{\kappa_B}$ Index Set/Leading Order \\
\hline
$F_{11100}$ &  $-\frac{1}{2} + B_{110}$ \\
$F_{01100}$ &  $ -\frac{1}{2} + B_{220}$  \\
$F_{11000},$ $F_{10100,}$ $F_{22000},$ $F_{20200},$ $F_{d22}$ & $\infty$ \\
$F_{22200},$ $F_{02200}$ & $-\frac{n+2}{2} + B_{220}$\\
$F_{d3},$ $F_{d02}$ & $-\frac{n+3}{2} - k_b$ \\
$F_{10000},$ $F_{01000},$  $F_{00100},$ $F_{00010},$ $F_{00011}$ &$ \infty$  \\
\hline
\end{tabular}    

Now, recalling the formula:

$$(\beta_c)_* ( \tilde{\kappa_A } \tilde{\kappa_B }  (\rho_{11100} \rho_{11000} \rho_{01100} \rho_{10100})^{\frac{1}{2}} (\rho_{22200} \rho_{22000} \rho_{02200}  )^{\frac{n+2}{2}} $$
$$(\rho_{20200}) ^{\frac{n}{2}} (\rho_{d3} \rho_{d20} \rho_{d02} )^{\frac{n+3}{2}} (\rho_{d22}) ^{\frac{n+1}{2}} \rho_{01000}  \rho_{00011} \rho_{00010} \rho_{00001} (^b \nu_3 ^2 ) )$$

We see that the quantity on the right hand side to be pushed forward by $(\beta_c )_*$ has the following indices on the boundary faces.

\begin{tabular}{|l|l|}
\hline
Face & Index Set/Leading Order \\
\hline
$F_{11100}$ &  $-\frac{1}{2} + A_{110} + B_{110}$ \\
$F_{11000},$ $F_{01100},$ $F_{10100},$ $F_{22000},$ $F_{02200},$ $F_{20200}$ & $\infty$ \\
$F_{22200} $ & $-\frac{n+2}{2} + A_{220} + B_{220}$\\
$F_{d3}$ & $-\frac{n+3}{2} - (k_a+ k_b)$ \\
$F_{d20},$ $F_{d02},$ $F_{d22}$ & $\infty$  \\
$F_{10000},$ $F_{01000},$ $F_{00100},$  $F_{00010},$ $F_{00001},$ $F_{00011}$ &$ \infty$  \\
\hline
\end{tabular}    

The push forward under $(\beta_c )^*$ sends the boundary faces of $\bar{Z}^3 _{sc,h}$ to $\bar{Z}^2 _{sc,h}$ as follows.  

\begin{tabular}{|l|l|}
\hline
$\bar{Z}^3 _h$ Face & Boundary face of $\bar{Z}^2 _{sc,h}$ or Interior \\
\hline
$F_{11100}$ & $F_{110}$ \\
$F_{10100}$ & $F_{110}$ \\
$F_{22200},$ $F_{20200}$ & $F_{220}$ \\
$F_{d3},$ $F_{d22}$ & $F_{d2}$ \\
$F_{10000} $ & $F_{100}$ \\
$F_{00100}$  & $F_{010}$ \\
$F_{00011}$ & $F_{001}$  \\
$F_{11000},$ $F_{01100}$ $F_{22000},$  $F_{02200},$ & Interior\\
$F_{d20},$ $F_{d02},$ $F_{01000},$ $F_{00010},$ $F_{00001}$ & Interior  \\
\hline
\end{tabular}

The quantity to be pushed forward is integrable with respect to $^b \nu^2 _3$ at the faces that are mapped to the interior, so we may apply the push forward theorem (see \cite{edge}) to arrive at the result of the composition rule.  The kernel, $\kappa_{B \circ A}$ will have the following polyhomogeneous index sets and leading orders on $\bar{Z}^2 _h$.
 
\begin{tabular}{|l|l|}
\hline
Face of $\bar{Z}^2 _h$ & Index Set/Leading Order \\
\hline
$F_{110}$ & $ - \frac{1}{2} + A_{110} + B_{110} $ \\
$F_{220}$ & $-\frac{n+2}{2} + A_{220} + B_{220} $ \\
$F_{d2}$ & $-\frac{n+3}{2} - (k_a + k_b) $\\
$F_{100}$ & $\infty$ \\
$F_{010}$ & $\infty$ \\
$F_{001}$ & $\infty $ \\
\hline
\end{tabular}

This concludes the proof of the composition rule:  $B \circ A$ is an element of $\Psi^{A_{110} + B_{110}, A_{220} + B_{220}, k_a + k_b} _{ac, H}$. 
\begin{flushright} $\heartsuit$ \end{flushright}

\begin{acparametrix}
Let $(Z, g_z)$ be an asymptotically conic scattering space with cross section $(Y, h)$ at infinity.  Let $(E, \nabla)$ be a Hermitian vector bundle over $(Z, g_z)$ so that near the boundary $E$ is the pullback of a bundle over $(Y, h).$  Let $\Delta$ be a geometric Laplacian on $(Z, g_z)$ associated to the bundle $(E, \nabla).$   Then there exists $H \in \Psi^{E_{110}, E_{220} , -2} _{ac, H} $ satisfying:  

$$ (\partial_t + \Delta) H(z, z', t) = 0,\,  t>0, $$
$$ H(z, z', 0) = \delta (z-z').$$
Moreover, $H$ vanishes to infinite order at $F_{110}$ and is smooth up to $F_{220}.$  
\end{acparametrix}  

On the interior of $\bar{Z}^2 _{ac, h}$ the ac scattering model heat kernel is locally defined by the Euclidean heat kernel and a partition of unity.  At $F_{d2}$ we construct the model heat kernel explicitly using the jet of the metric at the base point of each fiber.  At $F_{001}$ the model heat kernel vanishes to infinite order.  At $F_{110}$ and $F_{220}$ the model heat kernel is the lift of the Euclidean heat kernel.  Then the ac scattering model heat kernel $H_1$ satisfies
$$(\partial_t + \Delta) H_1 = K_1,$$
where $K_1$ vanishes to positive order on the boundary faces of $\bar{Z}^2 _{sc, h}.$  We now define
$$H_2 = H_1 - H_1  K_1,$$
with 
$$(\partial_t + \Delta) H_2 = K_2$$
where $K_2$ vanishes to one order higher on the boundary faces of $\bar{Z}^2 _{sc, h}.$  Similarly, 
$$H_3 := H_2 - H_2  K_2.$$
Using Borel summation we construct $H_{\infty}$ with expansion asymptotic to $H_1, H_2, H_3, \ldots$ and satisfying
$$(\partial_t + \Delta)H_{\infty} = K,$$
where now $K$ vanishes to infinite order on the boundary faces of $\bar{Z}^2 _{sc, h}.$   As a $t$-convolution operator we wish to have 
$$H_{\infty} = Id,$$
however, we currently have
$$H_{\infty} = Id + K,$$
but this is not a problem since $(Id + K)$ is invertible with inverse of the same form.  Then the ac scattering heat kernel
$$H = H_{\infty} (Id - K)^{-1}$$
is an element of the ac scattering heat calculus with leading orders on the boundary faces of $\bar{Z}^2 _{sc, h}$ determined by those of the model kernel. 
\begin{flushright} $\heartsuit$ \end{flushright}


\section{Acc triple heat space}

Let 
$$T  : = [ [ [ [ \cS \times \cS \times \cS ; Y \times Y' \times Y''] ; Y \times Y' ] ; Y' \times Y ''] ; Y \times Y'' ],$$
where we have used $Y, Y', Y''$ to denote the three copies of $Y$ in $\cS ^3.$  Let 
$$\cT := \{ p \in T : f_i (p) = 0, \, i=1,2 \}, \, f_1 (p) = x(p) r(p) - x'(p) r'(p), \, f_2 (p) = x'(p) r'(p) - x''(p) r''(p).$$
Like the acc double and heat space, $\cT$ is a smooth manifold with corners.  Let 
$$\R^+ _{2,b} := [ \R^+ _t \times \R^+ _s ; \{0\} \times \{0\} ].$$
Then the acc triple heat space is constructed from $\cT \times \R^+ _{2,b}$ by blowing up along twelve  submanifolds creating the following twelve boundary faces.  Below, let $tD$ be the lift to $\cT$ of the diagonal in $\cS \times \cS ' \times \cS '',$ let $D_{110}$ be the lift of the diagonal in $\cS \times \cS ',$ $D_{011}$ be the lift of the diagonal in $\cS' \times \cS ''$ and $D_{101}$ be the lift of the diagonal in $\cS \times \cS''.$

\begin{tabular}{|l|l|}
\hline
Submanifold blown up & Face created \\
\hline
$Y \times Y' \times Y'' \times \{0\} \times \{0\}, dt, ds$ & $S_{11122}$ \\
$ Y \times Y' \times \{t=0\}, dt $ & $S_{11020},$ \\
$Y' \times Y'' \times \{s=0\}, ds$ &$S_{01102}$ \\
$Y \times Y'' \times \{t=s=0 \}, ds, dt$ & $S_{10122}$ \\
$Y \times Y' \times Y'' $ & $S_{111}$\\
$Y \times Y'$ & $S_{110}$\\
$Y' \times Y''$ & $S_{011}$\\
$Y \times Y''$ & $S_{101}$\\
$tD \times \{t,s=0\}, ds, dt$ & $S_{td}$ \\
$D_{110} \times \{t=0\}, dt$ & $S_{d20}$ \\
$D_{011} \times \{s=0\}, ds$ &$S_{d02}$\\
$D_{101} \times \{s=t=0\}, ds, dt$ & $S_{d22}$\\

\hline
\end{tabular}

As constructed, the acc triple heat space has full and partial projection/blow down maps to three identical copies, left, right and center, of the acc heat space and to three corresponding copies of the blown down space $\{ \epsilon = \epsilon ' \} \subset \cS^2 \times \R^+.$  To compose two elements $A$ and $B$ we view the element $A$ as acting from the left to the right while $B$ acts from the right to the center.  Formally, the composition $B \circ A$ is the pushforward from the acc triple heat space of the product of the lifts of $A$ and $B.$  Compatibility assumptions on the leading orders of $A$ and $B$ at boundary faces of the acc heat space are required so that we can push forward.  With these assumptions and with the possible inclusion of normalizing factors at boundary faces of the acc heat space, two elements compose as one would expect.  The technical details in the proof of this composition rule are expected to be analogous to the technical details in the proof of the ac scattering heat calculus composition rule (appendix A).  

\end{appendix}


\begin{thebibliography}{99}

\bibitem{pierre}  P. Albin, {\em The Gauss-Bonnet theorem and Index theory on conformally compact manifolds.}  Ph.D. Dissertation, Stanford University, June 2005.

\bibitem{chavel} I. Chavel, \em Eigenvalues in Riemannian Geometry.  \em Academic Press, 1984. 

\bibitem{chee0} J. Cheeger, \em On the spectral geometry of spaces with cone-like singularities, \em Proc. Nat. Acad. Sci. U.S.A. 76 (1979), 2103-2106.

\bibitem{chee} J. Cheeger, \em Spectral Geometry of Singular Riemannian Spaces. \em  J. Diff. Geo., 18, (1984) 575-657.

\bibitem{cc1} J. Cheeger and T. Colding, \em On the structure of spaces with Ricci curvature bounded below, I.  \em J. Diff. Geo., 46 (1997), no.3, 406-480.

\bibitem{cc2} J. Cheeger and T. Colding, \em On the structure of spaces with Ricci curvature bounded below; II. \em J. Diff. Geo., 54 (2000), no. 1, 13-36.

\bibitem{cc3} J. Cheeger and T. Colding, \em On the structure of spaces with Ricci curvature bounded below. III. \em J. Diff. Geo.  54  (2000),  no. 1, 37-74

\bibitem{ding} Y. Ding, \em Heat kernels and Green's functions on limit spaces. \em Communications in Analysis and Geometry, 10, no. 3 (2002) 475-514.

\bibitem{epmm} C. Epstein, R. Melrose and G. Mendoza, \em The Heisenberg algebra, index theory and homology.  \em work-in-progress, http://math.mit.edu/~rbm/book.html.

\bibitem{fukaya}  K. Fukaya, \em Collapsing of Riemannian manifolds and eigenvalues of Laplace operator. \em  Invent. Math.  87, no. 3,  (1987),  517-547.

\bibitem{ga} M. Gaffney, \emph{A special Stokes' Theorem for complete Riemannian manifolds,} Ann. of Math., \textbf{60} (1954), 140-145.

\bibitem{gil} J. Gil,  \emph{Full Asymptotic Expansion of the Heat Trace for Non-Self-Adjoint Elliptic Cone Operators}, Math. Nachr. 250 (2003), 25--57.   

\bibitem{domains} J. Gil and G. Mendoza, \em Adjoints of elliptic cone operators, \em  Amer. J. Math. 125 (2003) 2, 357--408.  

\bibitem{hv} A. Hassell and A. Vasy,  \em The resolvent for Laplace-type operators on  asymptotically conic spaces.  \em Ann. l'Inst. Fourier 51 (2001), 1299-1346

\bibitem{hmm} A. Hassell, R. Mazzeo and R. Melrose,  \em A signature formula for manifolds with corners of codimension 2. \em  Topology, 36 No. 5 (1997), pp. 1055-1075.

\bibitem{j} D. Joyce,  \emph{Compact Manifolds with Special Holonomy}, Oxford University Press, 2000.  

\bibitem{k} S. Krantz,  \emph{Partial Differential Equations and Complex Analysis}, CRC Press, 1992.  

\bibitem{le} N. Lebedev,  \emph{Special functions and their applications}, Prentice-Hall, Inc., 1965.

\bibitem{lesch} M. Lesch, \em Differential operators of Fuchs type, conical singularities, and asymptotic methods, \em Teubner Texte zur Mathematik Vol. 136, Teubner--Verlag, Leipzig, 1997.

\bibitem{lott1} J. Lott, \em Remark about the spectrum of the p-form Laplacian under a collapse with curvature bounded below. \em Proc. of the AMS 132, p. 11-918, 2003.  

\bibitem{lott2} J. Lott, \em Collapsing and the Differential Form Laplacian : The Case of a Singular Limit Space. \em preprint, http://www.math.lsa.umich.edu/~lott.

\bibitem{lo} P. Loya:  \emph{Asymptotic Properties of the Heat Kernel on Conic Manifolds},  Israel J. Math.136  (2003), 285--306.

\bibitem{edge} R. Mazzeo: \emph{Elliptic Theory of Differential Edge Operators I},  Comm. Partial Differential Equations  16  (1991),  no. 10, 1615--1664.

\bibitem{eta} R. Mazzeo and R. Melrose, \em Analytic Surgery and the Eta Invariant. \em  Geom. Funct. Anal.  5  (1995),  no. 1, 14--75.

\bibitem{pat} P. T. McDonald, \em The Laplacian for Spaces with Cone-Like Singularities. \em Massachusetts Institute of Technology, 1990. 

\bibitem{tapsit}R. Melrose: \emph{The Atiyah-Patodi-Singer Index Theorem} Research Notes in Mathematics, 4. A K Peters, Ltd., 1993.

\bibitem{moo} E. Mooers: \emph{Heat Kernel Asymptotics on Manifolds with Conic Singularities},  J. Anal. Math.  78  (1999), 1--36.

\bibitem{rose} S. Rosenberg: \emph{The Laplacian on a Riemannian Manifold}, Cambridge University Press, 1997.

\bibitem{moi} J. M. Rowlett:  \emph{Spectral Geometry and Asymptotically Conic Convergence}, Ph.D. Dissertation, Stanford University, June, 2006.  

\bibitem{shu} M. A. Shubin:  \emph{Pseudodifferential Operators and Spectral Theory}, Springer, 2000.  


\end{thebibliography}
\end{document}